\documentclass[10pt]{amsart}

\allowdisplaybreaks

\usepackage[T1]{fontenc}
\usepackage[utf8]{inputenc}
\usepackage{amsmath,amssymb,amsfonts,mathtools}
\usepackage{enumitem}
\usepackage{hyperref}
\usepackage{xcolor}
\usepackage{multicol}
\usepackage[normalem]{ulem}

\usepackage{xltabular,array}
\newcounter{lowercriterion}
\newcounter{uppercriterion}
\renewcommand{\thelowercriterion}{(\roman{lowercriterion})}
\renewcommand{\theuppercriterion}{(\roman{uppercriterion}')}
\newcommand{\cellmath}[1]{%
  \par\smallskip
  \begingroup
    \centering
    \(\displaystyle #1\)%
    \par
  \endgroup
  \smallskip
}

\hypersetup{
  colorlinks=true,
  linkcolor=blue,
  citecolor=blue,
  urlcolor=blue
}

\newtheorem{theorem}{Theorem}[section]
\newtheorem{proposition}[theorem]{Proposition}
\newtheorem{lemma}[theorem]{Lemma}

\newtheorem{assumption}{Assumption}[section]

\theoremstyle{definition}

\newtheorem{example}{Example}[section]
\theoremstyle{remark}
\newtheorem{remark}{Remark}[section]

\DeclareMathOperator{\Var}{Var}

\newcommand{\F}{\mathcal F}
\newcommand{\E}{\mathbb E}
\newcommand{\Pbb}{\mathbb P}
\newcommand{\ind}{\mathbf 1}

\newcommand{\eps}{\varepsilon}

\newcommand{\htr}{r^*}

\newcommand{\Deltac}{\delta_c}

\title[Exact a.s. scale for non-negative recursion]{Exact 
almost sure scale for \\non-negative stochastic recursions}

\author[G. Aletti]{Giacomo Aletti}
\address{Universit\`a degli Studi di Milano, Milan, Italy}
\email{giacomo.aletti@unimi.it}

\author[I. Crimaldi]{Irene Crimaldi}
\address{IMT School for Advanced Studies Lucca, Lucca, Italy}
\email{irene.crimaldi@imtlucca.it}

\author[A. Ghiglietti]{Andrea Ghiglietti}
\address{Universit\`a degli Studi di Milano-Bicocca, Milan, Italy}
\email{andrea.ghiglietti@unimib.it (Corresponding author)}

\thanks{G.\ Aletti is a member of ``Gruppo Nazionale per il Calcolo
  Scientifico (GNCS)'' of the Italian Institute ``Istituto Nazionale
  di Alta Matematica (INdAM)''. I.\ Crimaldi and A.\ Ghiglietti are members of ``Gruppo
  Nazionale per l'Analisi Matematica, la Probabilit\`a e le loro
  Applicazioni (GNAMPA)'' of the Italian Institute ``Istituto
  Nazionale di Alta Matematica (INdAM)''. 
  I.\ Crimaldi and  A.\ Ghiglietti are partially supported by Prin 2026 –
DD n. 2298, 10-04-2026,  {\em Global Firms and Knowledge Networks under Geoeconomic Fragmentation (KNOWFRAG)}, funded by the Italian Ministry of University and Research.
  I.\ Crimaldi thanks the “Resilienza Economica e Digitale” (RED) 
project (CUP
D67G23000060001) funded by the Italian Ministry of University and Research (MUR)
as “Department of Excellence” (Dipartimenti di Eccellenza 2023-2027, Ministerial
Decree no. 230/2022).}


\begin{document}

\begin{abstract} 
We study exact almost sure scale for non-negative stochastic recursions whose conditional mean consists of a linear
term and a signed adapted perturbation. The linear term determines a
product normalization, but convergence after normalization does not
ensure a strictly positive limit. Multiplying the state by
the inverse of each deterministic contraction factor gives a
contraction-adjusted state. On paths where this state tends to infinity,
we prove that the logarithm of the normalized state differs from the
cumulative relative perturbation by an almost surely convergent
quantity. The assumptions control the effective step sizes, the quadratic positive relative
perturbation, and the conditional quadratic variation of the relative
martingale noise. If the positive perturbation is summable at the
normalized scale, the normalized process converges to a finite
non-negative limit; strict positivity is then equivalent to absolute
summability of the relative perturbation. 
We also give criteria for divergence of the contraction-adjusted state, first-order comparisons, and counterexamples
illustrating the role of the hypotheses. Harmonic step sizes yield
exact power laws with strictly positive random coefficients.
Applications include learning automata, urns with signed sampling bias,
response-adaptive reinforcement, generalized Simon models with signed
updates, branching processes in varying environments, and linear
stochastic approximation with multiplicative noise.
\end{abstract}

\maketitle

\section{Introduction}
Non-negative stochastic recursions arise in reinforcement models (urn models, innovations processes, etc.), 
learning algorithms, and population dynamics, which have applications in several fields, such as Economics and  Social Science, Physics and Computer Science, Biology and Medicine.  
In many such models, a linear conditional drift identifies a candidate
asymptotic scale. A remaining question is whether this scale is exact
along almost every sample path: after normalization, does the process
converge to a finite, strictly positive random variable? This question
is especially relevant where the state and the noise
may vanish together and convergence of the normalized process alone
does not exclude a zero limit. Strictly positive almost sure limits after normalization have been
proved in specific innovation, urn, and response-adaptive models
\cite{AlettiCrimaldiGhigliettiUnpublished,AlettiCrimaldiGhigliettiNutarelli2026,Janson2026,
MayFlournoy2009}.
For the two-armed bandit algorithm, such limits have also been
established on specified events and under suitable step-size
conditions \cite{LambertonPages2008, LambertonPagesTarres2004}; 
while, in a framework of interacting agents,  \cite{AlettiCrimaldiGhiglietti2017,AlettiCrimaldiGhiglietti2024Polarization} studies when the common limit agents' inclination takes 
 values in the open interval $(0,1)$ or/and in the extremes $\{0,1\}$.  
All these model-specific results motivate a general criterion for
strict positivity in the recursions considered here.\\ 
\indent We consider non-negative recursions of the form
\begin{equation}\label{eq:intro-recursion}
 X_0\geq 0,\qquad   X_{t+1}=(1-p_{t+1})X_t+r_{t+1}Y_{t+1},
\qquad
\E(Y_{t+1}\mid\F_t)=\delta^*X_t+\rho_t,
\end{equation}
where $Y_{t+1}\geq0$, $\delta^*>0$, and $\rho_t$ is an adapted,
possibly signed perturbation. The deterministic coefficients
\(p_t\in[0,1)\) and \(r_t\geq0\) describe contraction and reinforcement, respectively. 
We also allow for the inclusion of an additional term $r_{t+1}X_t\Delta N_{t+1}$ in~\eqref{eq:intro-recursion}, 
where $(N_t)$ is a martingale, provided that the resulting state remains
non-negative. This extension permits signed effective increments even
though the reinforcement variable $Y_{t+1}$ is non-negative.
\\
\indent The linear part of the conditional mean suggests the normalization
\[
q_t^{(\delta^*)}
:=\prod_{k=1}^{t}(1-p_k+\delta^*r_k)^{-1},
\qquad X_t^*:=q_t^{(\delta^*)}X_t\,,
\]
so that the linear drift $\delta^*X_t$ is cancelled and we get $\E(\Delta X^*_{t+1}\mid\F_t)=q_{t+1}^{(\delta^*)}r_{t+1}\rho_t$. 
Summability of the positive part of this quantity, i.e.
\begin{equation}\label{eq:summability_rho_+}
    \sum_tq_{t+1}^{(\delta^*)}r_{t+1}\rho_t^+<+\infty,
\end{equation}
therefore gives a real non-negative almost sure limit by the
Robbins-Siegmund theorem. This argument, however, leaves open whether
the limiting random variable is strictly positive. The distinction
separates an upper asymptotic bound from an exact asymptotic
equivalence.
\\
\indent Our main results resolve this question through a logarithmic
comparison. Set
\[
\pi_t:=\prod_{k=1}^{t}(1-p_k)^{-1},
\qquad X^{(\pi)}_t:=\pi_tX_t,
\qquad \htr_t:=\frac{r_t}{1-p_t}.
\]
For the basic recursion~\eqref{eq:intro-recursion}, $X^{(\pi)}_t$ is non-decreasing and so the sample space can be split into two alternatives: 
$\{\lim_t X^{(\pi)}_t<+\infty\}$ or $\{\lim_t X^{(\pi)}_t=+\infty\}$. Since on the first event we simply have $X_t\sim X^{(\pi)}_\infty/\pi_t$ 
(provided $(X_t)$ is not identically zero),
we will focus on the behavior of $(X_t)$ on the event $\{X^{(\pi)}_t\to+\infty\}$. 
With the additional martingale noise $r_{t+1}X_t\Delta N_{t+1}$, the monotonicity of $X^{(\pi)}_t$
is lost, so in general the previous alternatives no
 longer exhaust the sample space; however, we provide
sufficient criteria for both settings to cover all the sample space. Of course, it still holds true that on the event
$\{\limsup_t X^{(\pi)}_t < +\infty\}$ we have $X_t = O(1/\pi_t)$. Hence, our main analysis concerns the
asymptotic behavior of $(X_t)$ on $\{X^{(\pi)}_t \to +\infty\}$. It is also worthwhile to note that we will show some conditions, usually satisfied in applicative contexts,  
that imply this event to be almost sure.
\\
\indent The analysis on the event $\{X^{(\pi)}_t\to+\infty\}$ will be conducted under some usual hypothesis 
on the sequence $(\htr_t)$ of the effective step size
and a control on the predictable quadratic variation of the martingale
term involving both $Y_{t+1}$ and $\Delta N_{t+1}$.
Then, on $\{X^{(\pi)}_t\to+\infty\}$, a central role will be played by the cumulative relative perturbation
\[
\Gamma_T:=\sum_{t\leq T}
\htr_{t+1}\frac{\rho_t}{X_t},
\]
(with the sum starting from a sufficiently large random time (inclusive) after which the
state $X_t$ remains strictly positive). 
Under the hypothesis of square summability of the positive relative perturbation increments, 
that is  $\sum_{t}
(\htr_{t+1}{\rho_t^+}/{X_t})^2<+\infty$, we prove 
\[
\ln X_{T+1}^*-\Gamma_T\stackrel{a.s.}\longrightarrow L
\qquad\text{on }\{X^{(\pi)}_t\to+\infty\},
\]
where $L$ is a real random variable. Thus, on the considered event, the cumulative relative perturbation $\Gamma_T$
describes the asymptotic behavior of the process $(X^*_t)$, that is of the process $(X_t)$ after the scaling for removing  
the linear drift. In particular, it characterizes whether the normalized state
has subsequences approaching zero or diverging to infinity: specifically, it establishes if 
\[
0 < \liminf_{T\to\infty} X_T^*
\leq \limsup_{T\to\infty} X_T^* < +\infty.
\] 
In particular, for the basic recursion~\eqref{eq:intro-recursion}, a 
conditional variance bound on $Y_{t+1}$ ensures the
required control of the predictable quadratic variation of the
martingale term.
Hence, under this conditional variance bound and the effective-step assumptions, we show on $\{X^{(\pi)}_t\to+\infty\}$
that the total variation of the cumulative relative perturbation implies the exact rate: that is,
\[\sum_t\htr_{t+1}\frac{|\rho_t|}{X_t}<+\infty
\quad\Longrightarrow\quad X_\infty^*:=\lim_t X^*_t>0.
\]
(The two conditions become equivalent when~\eqref{eq:summability_rho_+} holds.)
\\
\indent For asymptotically harmonic coefficients,
$p_t=a/(t+\theta)+O(t^{-1-\eta})$ and
$r_t=b/(t+\theta)+O(t^{-1-\eta})$, with $\eta>0$, the product
$q_t^{(\delta^*)}$ is asymptotically equivalent to a strictly positive real constant times
$t^{a-b\delta^*}$.  
Hence, the same framework covers decay, convergence to a positive limit,
and growth. More general step sizes retain their exact product scale,
which need not be a power or even converge.
\\
\indent The section of applications illustrate both recovery of established rates and
extensions of a large number of familiar models. Among others, they include a linear reward-inaction
algorithm, a triangular urn with signed predictable sampling bias,
a randomly reinforced response-adaptive urn, and a generalized Simon
model whose scores can decrease. A branching-process example has a
strictly positive power-law limit on survival despite infinitely many
negative population increments. A final example treats a positive
linear recursion with multiplicative martingale noise. The common
argument identifies the scalar linear scale, verifies the a.s. divergence of $(X^{(\pi)}_t)$, and controls the relative
noise and perturbation series.

\subsection*{Structure}
The paper is organized as follows. Section~\ref{sec:model} 
introduces the basic stochastic recursion, its extension with multiplicative martingale noise, and the relevant assumptions. Section~\ref{sec:main-last} 
presents the main logarithmic comparison, the criterion for strict positivity of the normalized limit, and consequences useful for applications. 
Section~\ref{sec:proofs} contains the proofs. Section~\ref{sec:lower} develops complementary first-order criteria for the basic recursion, while Section~\ref{sec:rates} 
provides examples and counterexamples illustrating the possible asymptotic regimes and the role of the hypotheses. Section~\ref{sec:applications} discusses applications to learning automata, reinforced urns, generalized Simon models, branching processes in varying environments, and linear stochastic approximation with multiplicative noise. Finally, Appendix~\ref{sec:upper} 
establishes the existence of a finite normalized limit, and Appendix~\ref{sec:technical} contains a technical lemma.

\section{Framework}\label{sec:model}

This section provides the basic stochastic recursion of interest and its ge\-ne\-ra\-li\-zation with a multiplicative martingale noise. 
The used assumptions are stated according to their role.\\

\indent Throughout, \((\Omega,\F,(\F_t)_{t\ge0},\Pbb)\) is a discrete-time
filtered probability space and we set $a^+:=\max(a,0)$ and $a^-:=\max(-a,0)$. Empty sums are zero and empty products are one.  The notation 
$A_t=O(B_t)$ means that, with probability one, $B_t$ is eventually non-negative and $|A_t|\leq CB_t$ eventually, for a suitable non-negative real random variable $C$; 
$A_t=o(B_t)$ and $A_t\stackrel{a.s.}\sim B_t$ mean, respectively, $A_t/B_t\stackrel{a.s.}\to 0$ and $A_t/B_t\stackrel{a.s.}\to 1$ 
on the event where the ratio is eventually defined; and $A_t\asymp B_t$ means both
$A_t=O(B_t)$ and $B_t=O(A_t)$.

\subsection{Basic stochastic recursion} We focus on the non-negative recursion
\begin{equation}\label{eq:recursion}
        X_0\ge0,\qquad
        X_{t+1}=(1-p_{t+1})X_t+r_{t+1}Y_{t+1},
        \qquad t\ge0,
\end{equation}
where $X_0$ is $\F_0$-measurable, $Y_{t+1}\ge0$ is
$\F_{t+1}$-measurable, \(p_t\in[0,1)\) and \(r_t\geq0\) are deterministic and 
describe contraction and reinforcement, respectively. 

We set $\pi_0:=1$ and for $T\geq 0$
\begin{equation}\label{eq:P-H-V}
        \pi_{T+1}:=\prod_{t=1}^{T+1}\frac{1}{1-p_{t}}>0,\qquad
        X^{(\pi)}_T:=\pi_TX_T\geq 0\,.
\end{equation}
Moreover, we define for $t\geq 0$
\begin{equation}\label{eq:V}
        V_{t+1}:=\pi_{t+1}r_{t+1}Y_{t+1}\geq 0.
\end{equation}
Then, from \eqref{eq:recursion} we obtain 
\begin{equation}\label{eq:H-increment}
        X^{(\pi)}_{t+1}=X^{(\pi)}_t+V_{t+1},\quad\mbox{that is }
        X^{(\pi)}_{T+1}=X^{(\pi)}_0+\sum_{t\leq T}V_{t+1}.
\end{equation}
In particular, $(X^{(\pi)}_t)$ is a non-decreasing non-negative 
adapted process 
and so its almost sure limit $X^{(\pi)}_\infty=\lim_t X^{(\pi)}_t$ always exists and coincides with $\sup_t X^{(\pi)}_t$.
Therefore the events
\[
\{X^{(\pi)}_\infty<+\infty\}
\qquad\text{and}\qquad
\{X^{(\pi)}_t\to+\infty\}
\]
form an exhaustive dichotomy for the basic recursion.

\begin{remark}[The behavior on $\{X^{(\pi)}_\infty<+\infty\}$]\label{rem:bounded-H-branch} 
Note that, on the event $\{X^{(\pi)}_\infty=\lim_t X^{(\pi)}_t<+\infty\}$ we have $X_t=O(1/\pi_t)$. 
As a consequence, 
since $x \leq -\ln(1-x)$ for $0\leq x<1$, under the standard assumption 
$\sum_t p_t=+\infty$, we have $\pi_t\to +\infty$ and so 
$\{X^{(\pi)}_\infty=\lim_t X^{(\pi)}_t<+\infty\}\subseteq \{X_t\to 0\}$.
Moreover, if $X_0>0$, then $X^{(\pi)}_0>0$ and so $X^{(\pi)}_\infty>0$. This means that, 
on the event $\{X^{(\pi)}_\infty=\lim_t X^{(\pi)}_t<+\infty\}\cap\{X_0>0\}$ we have $X_t\stackrel{a.s.}\sim X^{(\pi)}_\infty/\pi_t$. 
In particular, when $\sum_t p_t<+\infty$, we have that on the latter event $X_t$ almost surely converges to a strictly positive real random variable; while, 
if $\sum_t p_t=+\infty$, then  $X_t$ converges almost surely toward zero with exact velocity~$1/\pi_t$.
\end{remark}

Given the above remark, the main case of interest is the asymptotic behavior of the process $(X_t)_t$ on $\{X^{(\pi)}_t\to +\infty\}$. 
In the sequel, on the event \(\{X^{(\pi)}_t\to+\infty\}\), sums containing
\(1/X^{(\pi)}_t\) or \(1/X_t\) start from the first strictly positive term;
this changes only finitely many summands.

For $\delta> 0$, we set $q_0^{(\delta)}:=1$ and for $T\geq 0$
\begin{equation}\label{eq:Q}
        q^{(\delta)}_{T+1}:=
        \prod_{t=1}^{T+1}\frac{1}{1-p_{t}+\delta r_{t}}.
\end{equation}

\begin{remark}[When the deterministic normalization diverges]\label{rem:Q-infty} Since $x \leq -\ln(1-x)$ for $0\leq x<1$, if 
\begin{equation}\label{eq-Q-inf}
p_t-\delta r_t \geq 0\;\mbox{eventually and }
\sum_t (p_t-\delta r_t)=+\infty, 
\end{equation} 
then  we have $q_t^{(\delta)}\to +\infty$. Note that this is not an assumption of this paper. We include also cases when, for the value $\delta^*$ of interest we have 
$q_t^{(\delta^*)}\not\rightarrow +\infty$, 
see Section~\ref{subsec:Q-regimes} for a discussion.
\end{remark}

We now introduce the assumptions that we are going to use in the results for the basic dynamics. 
We group them according to their role.  The first assumption
introduces the drift and the perturbation terms.  The second concerns the effective time scale.
The third controls the conditional fluctuations.  The fourth is the
second-order condition needed for an additive logarithmic comparison,
whereas the fifth is used only when convergence of the normalized
process is required.

\begin{assumption}[Linear conditional mean with perturbation]\label{ass:reinforcement}
There exists $\delta^*>0$ and an adapted real-valued process $(\rho_t)$
such that, for every $t\geq0$,
\begin{equation}\label{eq:condmean}
        \E(Y_{t+1}\mid\F_t)=\delta^*X_t+\rho_t\qquad\text{a.s.}\ .
\end{equation}
The term $\delta^* X_t$ is a linear drift, while $\rho_t$ represents a perturbation term.
In particular, since $Y_{t+1}\geq 0$, we necessarily have \(\rho_t^-\leq\delta^*X_t\) and 
\begin{equation}\label{eq:bound-rho-meno}
\frac{\rho_t^-}{X_t}\leq\delta^*\qquad\mbox{on } \{X_t>0\}\,.
\end{equation}
\end{assumption}
Under Assumption~\ref{ass:reinforcement}, write
\[
X_t^*:=q_t^{(\delta^*)}X_t.
\]
Define the effective step size and cumulative effective time by
\[
    \htr_t:=\frac{r_t}{1-p_t},
    \qquad
    s_{T+1}:=\sum_{t\leq T}\htr_{t+1}.
\]

\begin{assumption}[Effective step size]\label{ass:step-size}
The sequence $(\htr_t)$ satisfies
\begin{enumerate}[label=\textup{\ref{ass:step-size}\alph*)},
                  ref={\ref{ass:step-size}\alph*}]
\item\label{ass:step-size-1a}
$\htr_t\to 0$;
\item\label{ass:step-size-1b}
$s_T=\sum_{t< T}\htr_{t+1}\to +\infty$;
\item\label{ass:step-size-2}
$\sum_t(\htr_t)^2<+\infty$.
\end{enumerate}
\end{assumption}
{Since \(0\leq r_t=(1-p_t)\htr_t\leq\htr_t\),
Assumptions~\ref{ass:step-size-1a} and
\ref{ass:step-size-2} imply, respectively, \(r_t\to0\)
and \(\sum_t r_t^2<+\infty\). 
\begin{assumption}[Conditional-variance bound]\label{ass:moment}
The increments in \eqref{eq:V} satisfy
\begin{equation}\label{eq:moment}
\Var(V_{t+1}\mid\F_t)
=O\bigl(\E(V_{t+1}\mid\F_t)\bigr),
\end{equation}
or, equivalently,
\[
\pi_{t+1}r_{t+1}\Var(Y_{t+1}\mid\F_t)
=O\bigl(\E(Y_{t+1}\mid\F_t)\bigr).
\]
\end{assumption}

\begin{remark}[A convenient conditional-variance criterion]
\label{rem:Pr-condition}
Assumption~\ref{ass:moment} follows if
\begin{equation}\label{eq:Y2-bound}
\Var(Y_{t+1}\mid\F_t)
=O\bigl(\E(Y_{t+1}\mid\F_t)\bigr)
\end{equation}
and
\begin{equation}\label{eq:Pr-bounded}
\sup_t\pi_{t+1}r_{t+1}<+\infty.
\end{equation}
The latter condition is equivalent to
\[
r_{t+1}=O(1/\pi_{t+1})=O\Bigl(\prod_{s=1}^{t+1}(1-p_s)\Bigr)
\]
If $p_t\to0$ and $\sum_tp_t^2<+\infty$, it can also be written as
\[
r_{t+1}=O\Bigl(\exp \Bigl\{-\sum_{s=1}^{t+1}p_s\Bigr\}\Bigr).
\]
In particular, if $p_t=a/t+O(t^{-1-\eta})$, then
$\pi_t\asymp t^a$ and \eqref{eq:Pr-bounded} becomes
$r_t=O(t^{-a})$.  Thus, when $r_t\sim b/t$, it holds exactly for
$a\leq1$.
\end{remark}

\begin{remark}[Bounded reinforcement]\label{ex:bounded}
If $0\leq Y_{t+1}\leq \kappa$ a.s. for a deterministic real constant $\kappa$, then
\[
\Var(Y_{t+1}\mid\F_t)
\leq\E(Y_{t+1}^2\mid\F_t)
\leq \kappa\ \E(Y_{t+1}\mid\F_t).
\]
Hence \eqref{eq:Y2-bound} holds, and
Assumption~\ref{ass:moment} follows from \eqref{eq:Pr-bounded}.
This applies in particular to the indicator reinforcements in
Section~\ref{sec:applications}.
\end{remark}

\begin{assumption}[Quadratic relative-perturbation control]
\label{ass:quadratic-relative-perturbation}
\begin{equation}\label{eq:log_inc_mean_2}
\sum_t\Bigl(
\htr_{t+1}\frac{\rho_t^+}{X_t}
\Bigr)^2<+\infty
\qquad\text{a.s.\ on } \{X^{(\pi)}_t\to+\infty\}.
\end{equation}
\end{assumption}
\begin{remark} Note, under Assumptions~\ref{ass:reinforcement} and~\ref{ass:step-size-2}, Assumption~\ref{ass:quadratic-relative-perturbation} is equivalent to 
$\sum_{t}\left(\htr_{t+1}{\rho_t}/{X_t}\right)^2<+\infty$. Indeed, we have 
$$\begin{aligned}
\sum_{t\leq T}\left(\htr_{t+1}\frac{\rho_t^+}{X_t}\right)^2\leq \sum_{t\leq T}\left(\htr_{t+1}\frac{\rho_t}{X_t}\right)^2 &= \sum_{t\leq T}\left(\htr_{t+1}\frac{\rho_t^+}{X_t}\right)^2 + \sum_{t\leq T}\left(\htr_{t+1}\frac{\rho_t^-}{X_t}\right)^2\\
&\leq \sum_{t\leq T}\left(\htr_{t+1}\frac{\rho_t^+}{X_t}\right)^2 + (\delta^{*})^2\sum_{t\leq T}(\htr_{t+1})^2\,.
\end{aligned}$$
\end{remark}
\begin{assumption}[Summable normalized positive perturbation]
\label{ass:lim_sup2}
\begin{equation}\label{eq:rho-positive}
\sum_t q_{t+1}^{(\delta^*)}r_{t+1}\rho_t^+<+\infty
\qquad\text{a.s.}
\end{equation}
\end{assumption}
The above assumption implies $q_{t+1}^{(\delta^*)}r_{t+1}\rho_t^+\stackrel{a.s.}\to 0$. Hence, it is natural to ask if the positive perturbation  
$\rho_t^+$ itself is a.s. vanishing. The following remark shows that Assumption~\ref{ass:lim_sup2}, together with other suitable assumptions, implies $\liminf_t\rho_t^+=0$ a.s.\,.
\begin{remark}[A consequence of the summability of the normalized positive perturbation]
 If $\sum_t q_{t+1}^{(\delta^*)}r_{t+1}=+\infty$ and \eqref{eq:rho-positive} holds, then 
 we necessarily have $\liminf_t\rho_t^+=0$ a.s.\,.
Moreover notice that, since \[
q_{t+1}^{(\delta^*)}r_{t+1}
=q_t^{(\delta^*)}
\frac{\htr_{t+1}}{1+\delta^*\htr_{t+1}},
\]
the divergence of the series $\sum_t q_{t+1}^{(\delta^*)}r_{t+1}$ follows from 
Assumptions~\ref{ass:step-size-1a} and \ref{ass:step-size-1b}, provided 
$q_t^{(\delta^*)}$ is bounded away from zero (see, for instance, the case of Remark~\ref{rem:Q-infty}). 
\end{remark}

\subsection{Stochastic recursion with multiplicative martingale noise}

For the non-negative process $(X_t)$, we replace the dynamics~\eqref{eq:recursion} by the following modification: 
\begin{equation}\label{eq:recursion:mart}
X_0\geq 0,\quad
X_{t+1}
=(1-p_{t+1})X_t+r_{t+1} (Y_{t+1}+X_t\Delta N_{t+1}),
\end{equation}
where $X_0$ is a $\F_0$-measurable random variable, $Y_{t+1}\ge 0$ is
$\F_{t+1}$-measurable, \(p_t\in[0,1)\) and \(r_t\geq0\) are deterministic, and  
$(N_t)$ is an $(\F_t)$-martingale with increments 
$\Delta N_{t+1}:=N_{t+1}-N_t$. The sequences $\pi_t$ and $q_t^{(\delta)}$ are defined as above and 
we still set $X^{(\pi)}_t=\pi_t X_t$ and $V_{t+1}=\pi_{t+1}r_{t+1}Y_{t+1}$.
Then, in this extended setting,  we have 
\begin{equation}\label{eq:dynamicsH-mart}
X^{(\pi)}_{t+1}-X^{(\pi)}_t
=\widetilde{V}_{t+1}:=V_{t+1}+\htr_{t+1}X^{(\pi)}_t\Delta N_{t+1},
\end{equation}
and, while $V_{t+1}\geq 0$, the process $\widetilde{V}_{t+1}$ is not necessarily non-negative. Hence $(X^{(\pi)}_t)$ 
is not necessarily non-decreasing. As a consequence, the previous alternatives ``bounded'' (i.e. $\lim_t X^{(\pi)}_t<+\infty$) and ``tends to $+\infty$'' (i.e. $X^{(\pi)}_t\to +\infty$) no longer
exhaust the sample space. Of course, it still holds true that on the event $\{\limsup_t X^{(\pi)}_t<+\infty\}$ we have $X_t=O(1/\pi_t)$. 
The following results concern the asymptotic behavior of \(X_t\)
on \(\{X^{(\pi)}_t\to+\infty\}\).

In the sequel, on the event $\{X^{(\pi)}_t\to+\infty\}$, all the quantities and the sums containing $1/X^{(\pi)}_t$ (and $1/X_t$) are understood starting from a sufficiently large random time (inclusive) after which we have $X^{(\pi)}_t>0$ (and so $X_t>0$); 
this changes only finitely many terms. 

\begin{assumption}[Summable relative martingale noise]
\label{ass:relative-noise:mart2}
\begin{equation}\label{eq:relative-bracket:mart2}
\sum_t\E\bigl((\Delta\mathcal M_{t+1})^2\mid\F_t\bigr)
<+\infty\qquad\mbox{a.s. on } \{X^{(\pi)}_t\to+\infty\}\,,
\end{equation}
where 
\begin{equation}\label{eq:relative-martingale-increment:mart2}
\Delta\mathcal M_{t+1}
:= \frac{\htr_{t+1}}{1+\delta^*\htr_{t+1}} 
\Bigl(\frac{\Delta M^Y_{t+1}}{X_t}+\Delta N_{t+1}\Bigr)
=
\frac{r_{t+1}}{a_{t+1}}\Bigl(
\frac{\Delta M^Y_{t+1}}{X_t}+\Delta N_{t+1}
\Bigr) 
\end{equation}
with  \(\Delta M^Y_{t+1}:=Y_{t+1}-\E(Y_{t+1}\mid\F_t)\) and 
\begin{equation}\label{eq:linear-coefficient:mart2}
a_{t+1}
:=1-p_{t+1}+\delta^*r_{t+1}
=(1-p_{t+1})(1+\delta^*\htr_{t+1}).
\end{equation}
\end{assumption}

\begin{remark}[A sufficient form of Assumption~\ref{ass:relative-noise:mart2}]
\label{rem:relative-noise:mart2}
By the conditional version of
\((x+y)^2\leq2x^2+2y^2\), Assumption
\ref{ass:relative-noise:mart2} follows from
\begin{equation}\label{eq:reinforcement-bracket:mart2}
\sum_t
\frac{r_{t+1}^2}
{a_{t+1}^2X_t^2}
\Var(Y_{t+1}\mid\F_t)
<+\infty\quad\text{a.s. on } \{X^{(\pi)}_t\to +\infty\}
\end{equation}
and
\begin{equation}\label{eq:N-relative-bracket:mart2}
\sum_t
\frac{r_{t+1}^2}{a_{t+1}^2}
\E((\Delta N_{t+1})^2\mid\F_t)
<+\infty\quad\text{a.s. on } \{X^{(\pi)}_t\to +\infty\}.
\end{equation}
Since
\(
    {r_{t+1}}/{a_{t+1}}
    =
    {\htr_{t+1}}/
    {(1+\delta^*\htr_{t+1})}
\),
under Assumption~\ref{ass:step-size-1a}, condition
\eqref{eq:N-relative-bracket:mart2} is equivalent to
\[
    \sum_t
    (\htr_{t+1})^2
    \E\bigl((\Delta N_{t+1})^2\mid\F_t\bigr)
    <+\infty\quad\text{a.s. on } \{X^{(\pi)}_t\to +\infty\}.
\]
Thus the relevant noise condition is weighted by the effective step
size. In particular, it follows from Assumption~\ref{ass:step-size-2}
whenever
\(\sup_t\E((\Delta N_{t+1})^2\mid\F_t)<+\infty\) a.s.\\
\indent Similarly, under Assumption~\ref{ass:step-size-1a}, condition~\eqref{eq:reinforcement-bracket:mart2} is equivalent to 
\begin{equation*}
\sum_t
(\htr_{t+1})^2\frac{\Var(Y_{t+1}\mid\F_t)}{X_t^2}
<+\infty\quad\text{a.s. on } \{X^{(\pi)}_t\to +\infty\}.
\end{equation*}
In particular, it follows from Assumption~\ref{ass:step-size-2}
whenever $\Var(Y_{t+1}\mid\F_t)=O(X_t^2)$.
\end{remark}

\begin{remark}[Direct reduction to the basic recursion]
\label{rem:effective-reinforcement:mart} Set $\widetilde{ Y}_{t+1}:=Y_{t+1}+X_t\Delta N_{t+1}$. 
Then \eqref{eq:recursion:mart} is exactly the basic recursion
\eqref{eq:recursion} with $Y_{t+1}$ replaced by
$\widetilde{Y}_{t+1}$, provided $\widetilde{Y}_{t+1}\geq 0$. Moreover, 
$\widetilde{V}_{t+1}=\pi_{t+1}r_{t+1}\widetilde{Y}_{t+1}$ and, 
since $(N_t)$ is a martingale, we have 
\[
\E(\widetilde{Y}_{t+1}\mid\F_t)=\delta^*X_t+\rho_t,
\qquad
\E(\widetilde{V}_{t+1}\mid\F_t)=\E(V_{t+1}\mid\F_t).
\]
Therefore, a convenient sufficient condition for the conditional-variance bound 
$$
\Var(\widetilde V_{t+1}\mid\F_t)=O\bigl(\E(\widetilde{V}_{t+1}\mid\F_t)\bigr).
$$
is
Assumption~\ref{ass:moment} together with
\begin{equation}\label{eq:relative-variance:mart}
(\htr_{t+1})^2 \bigl(X^{(\pi)}_t\bigr)^2
\E\bigl((\Delta N_{t+1})^2\mid\F_t\bigr)
=O\bigl(\E(V_{t+1}\mid\F_t)\bigr).
\end{equation}
Indeed, we have 
\begin{align*}
\Var(\widetilde V_{t+1}\mid\F_t)
&\leq2\Var(V_{t+1}\mid\F_t)
+2(\htr_{t+1})^2 \bigl(X^{(\pi)}_t\bigr)^2
\E\bigl((\Delta N_{t+1})^2\mid\F_t\bigr)\\
&=O\bigl(\E(V_{t+1}\mid\F_t)\bigr)
=O\bigl(\E(\widetilde{V}_{t+1}\mid\F_t)\bigr).
\end{align*}
\end{remark}

\section{Main results}\label{sec:main-last}
 Recall that $X^{(\pi)}_t=\pi_t X_t$ and 
\(X_t^*=q_t^{(\delta^*)}X_t\).  
 Set
\(
\tau:=\min\{n\geq 0: X_t>0\text{ for every }t\geq n\}.
\)
This is a random time (not necessarily a
stopping time), which, on the event \(\{X^{(\pi)}_t\to+\infty\}\), takes values in $\mathbb N$. 
 Throughout the paper, when working on \(\{X^{(\pi)}_t\to+\infty\}\), all the quantities and the sums involving \(1/X_t\) or
\(1/X^{(\pi)}_t\) start from \(\tau\); we write \(\sum_{t\le T}\)
for \(\sum_{t=\tau}^{T}\).}
On the event $\{X^{(\pi)}_t\to +\infty\}$, the effect of the perturbation $(\rho_t)$ is described on the logarithmic scale of $X_t$ by 
\begin{equation}\label{eq:Gamma}
\Gamma_T:=\sum_{t\leq T}\htr_{t+1}\frac{\rho_t}{X_t}=\sum_{t\leq T}\frac{\pi_{t+1}r_{t+1}\rho_t}{X^{(\pi)}_t}\,.
\end{equation}
Indeed, the following central result is an additive logarithmic comparison for the
normalized state $X^*_T$ with $\Gamma_T$; a second theorem turns convergence into a sharp
criterion for strict positivity of the limit. 

\begin{theorem}[Exact logarithmic comparison and two-sided control]
\label{thm:main-positive-all-maringale-noise}
Assume the modified recursion \eqref{eq:recursion:mart}.
Let Assumptions~\ref{ass:reinforcement}, \ref{ass:step-size},
\ref{ass:quadratic-relative-perturbation}, and \ref{ass:relative-noise:mart2} hold. Then,
on the event $\{X^{(\pi)}_t\to+\infty\}$, there exists a real random
variable $L$ such that
\begin{equation}\label{eq:logX-Gamma-additive}
    \ln X_{T+1}^*-\Gamma_T
    \xrightarrow{\mathrm{a.s.}} L .
\end{equation}
Then, on $\{X^{(\pi)}_t\to+\infty\}$, with probability one,
\[
    \liminf_T X_T^*=0
    \quad\Longleftrightarrow\quad
    \liminf_T\Gamma_T=-\infty,
\]
and
\[
    \limsup_T X_T^*=+\infty
    \quad\Longleftrightarrow\quad
    \limsup_T\Gamma_T=+\infty.
\]
\end{theorem}

The second result adds Assumption~\ref{ass:lim_sup2}.  The normalized
process then converges by Theorem~\ref{thm:finite-limit}, and the issue
is whether its limit is strictly positive.

\begin{theorem}[Strict positivity of the normalized limit]\label{thm:positive-all_second-maringale-noise}
Assume the modified recursion \eqref{eq:recursion:mart}.
Let Assumptions~\ref{ass:reinforcement}, \ref{ass:step-size}, \ref{ass:lim_sup2} and \ref{ass:relative-noise:mart2} hold. Then $(X^*_t)_t$ 
converges a.s.\, toward a non-negative real random variable $X^*_\infty$ and the following statements are equivalent on the event $\{X^{(\pi)}_t\to+\infty\}$:
\begin{enumerate}[label=\textup{(\Roman*)}]
    \item\label{pos-I2} $X^*_\infty>0$ a.s.;
    \item\label{pos-V2} $\liminf_T\Gamma_T>-\infty$ and
    Assumption~\ref{ass:quadratic-relative-perturbation} holds;
    \item\label{pos-VI2}
    \[
        \sum_{t}\htr_{t+1}\frac{|\rho_t|}{X_t}<+\infty.
    \]
\end{enumerate}
\end{theorem}

For the basic dynamics, the conditional-variance bound expressed in Assumption~\ref{ass:moment} implies the
relative-noise hypothesis required above as Assumption~\ref{ass:relative-noise:mart2}.

\begin{theorem}[The case of the basic recursion]\label{rem:conditional-variance-implies-relative-noise:mart2}
\label{thm:positive-all_second}\label{thm:main-positive-all}
Assume the basic recursion~\eqref{eq:recursion} (i.e. \eqref{eq:recursion:mart} with $N_t\equiv 0$).
Let Assumptions~\ref{ass:reinforcement}, \ref{ass:step-size} and \ref{ass:quadratic-relative-perturbation} hold.
Then, Assumption~\ref{ass:moment} implies
Assumption~\ref{ass:relative-noise:mart2}. Therefore, for the basic recursion, 
the above Theorems~\ref{thm:main-positive-all-maringale-noise} and~\ref{thm:positive-all_second-maringale-noise} 
hold true as well if we replace Assumption~\ref{ass:relative-noise:mart2} by Assumption~\ref{ass:moment}. 
\end{theorem}

\begin{remark}[Random adapted learning rates]
The deterministic character of $(p_t)$ and $(r_t)$ is not essential: indeed,  if 
$p_{t+1}$ and $r_{t+1}$ are $\F$-predictable and take
values in $[0,1)$, then $\pi_t$, $q_t^{(\delta)}$, and
$\htr_t=r_t/(1-p_t)$ are predictable processes as well, the definitions of 
$X^{(\pi)}_t$, $V_{t+1}$, and $\Gamma_T$ remain unchanged, and everything 
works well.
\end{remark}
\begin{remark}[Vector extensions and one-dimensional reductions] 
A fully vector-valued formulation is necessarily model dependent: positivity
does not select a unique scale, and one generally needs spectral, triangular, or
stability information to identify the leading coordinate.  A robust strategy is
therefore two-step.  First, vector theory identifies a positive one-dimensional
projection and proves that the transverse component is negligible at its scale.
Second, the projected coordinate is reduced to \eqref{eq:recursion} or~\eqref{eq:recursion:mart}, after which
the present scalar theorem supplies strict positivity of the normalized
coefficient, under suitable assumptions.  
\end{remark}

\subsection{Results useful for applications}\label{sec:byproducts}
In Theorem~\ref{thm:positive-all_second-maringale-noise}, Assumption~\ref{ass:lim_sup2} is not needed for the implication 
 \ref{pos-VI2}$\Rightarrow$\ref{pos-I2}, as shown in the following result.

\begin{theorem}[Absolute relative-perturbation criterion]
\label{prop:absolute-relative-criterion}
Assume the dynamics~\eqref{eq:recursion:mart}, together with
Assumptions~\ref{ass:reinforcement},
\ref{ass:step-size}, and
\ref{ass:relative-noise:mart2}. Then, on
\[
G:=\{X^{(\pi)}_t\to+\infty\}
\cap
\Bigl\{
 \sum_t
 \htr_{t+1}\frac{|\rho_t|}{X_t}<+\infty
\Bigr\},
\]
we have
\[
q_t^{(\delta^*)}X_t
\stackrel{a.s.}\longrightarrow
X_\infty^*\,,
\]
where $X_\infty^*$ is a strictly positive real random variable $X_\infty^*$. 
For the basic recursion~\eqref{eq:recursion}, the same conclusion holds
if Assumption~\ref{ass:relative-noise:mart2} is replaced by
Assumption~\ref{ass:moment}.
\end{theorem}

\begin{remark}[Asymptotically harmonic step-sizes]\label{cor:polynomial} The above results become particularly explicit when the step sizes are asymptotically harmonic. 
Indeed, the normalization $q_t^{(\delta^*)}$ can be turned into a power of $t$. 
Suppose that, for some $a\geq 0$, $b,\theta>0$ and $\eta>0$, we have 
\begin{equation}\label{eq:poly-rates}
        p_t=\frac{a}{t+\theta}+O(t^{-1-\eta}),\qquad
        r_t=\frac{b}{t+\theta}+O(t^{-1-\eta}).
\end{equation}
Then we can write 
\[
\begin{aligned}
        -\ln(1-p_t+\delta^*r_t)
        =p_t-\delta^*r_t+O(p_t^2+r_t^2)=\frac{a-b\delta^*}{t+\theta}+O(t^{-1-\eta'}),
\end{aligned}
\]
for some $\eta'>0$.  Since $\sum_{s\le t}{1}/{(s+\theta)}=\ln t+O(1)$, there exists a constant $\ell\in(0,+\infty)$ such that
$
        \frac{q_t^{(\delta^*)}}{t^{a-b\delta^*}}
        \longrightarrow \ell.
$
Therefore, $X^*_t\stackrel{a.s.}\longrightarrow X^*_\infty\in (0,+\infty)$ is equivalent to 
 $$
 t^{a-b\delta^*}X_t
        \stackrel{a.s.}\longrightarrow
        \frac{X^*_\infty}{\ell}\in (0,+\infty).
        $$
Note that the above a.s. convergence covers different asymptotic behaviors for $(X_t)$: 
decay to zero (when $a>b\delta^*$), convergence to a strictly positive real random variable (when $a=b\delta^*$), and growth to infinity (when $a<b\delta^*$).
\end{remark}

All the previous theorems are stated on the event $\{X^{(\pi)}_t\to+\infty\}$. Hence, the next result provides a very useful criterion to ensure that this event has probability one.  
Indeed, in particular,  in the case of the basic recursion it allows to conclude immediately that the event $\{X^{(\pi)}_t\to +\infty\}$ has probability one 
when the initial condition $X_0$ is a.s. strictly positive and the perturbation term $\rho_t$ is a.s. non-negative. This is a situation often met in applications 
(see, for instance, Section~\ref{sec:applications}).  The result below also deals with the more general recursion~\eqref{eq:recursion:mart}.

\begin{theorem}[Asymptotic behavior of \(X^{(\pi)}_t\)]
\label{cor-X0}\label{cor-X0:mart}\label{cor-X0:mart-signed}
\leavevmode
\begin{enumerate}[label=(\alph*)]
\item\label{cor-X0:a}
Consider the basic dynamics~\eqref{eq:recursion} under
Assumptions~\ref{ass:reinforcement} and~\ref{ass:moment}. Then \footnote{The sums are understood from the first time at which \(X_t>0\).}
\begin{equation}\label{eq-criterio-Hinfty2}
E_a:=\bigcup_t\{X_t>0\}
\cap
\left\{
\sum_t\htr_{t+1}
\left(\delta^*+\frac{\rho_t}{X_t}\right)=+\infty
\right\}
\subseteq
\{X^{(\pi)}_t\to+\infty\}
\quad\text{a.s.}
\end{equation}

If Assumption~\ref{ass:step-size-1b} holds as well, then
\[
\left\{\liminf_t\frac{\rho_t^+}{X_t}>0\right\}
\cup
\left\{\limsup_t\frac{\rho_t^-}{X_t}<\delta^*\right\}
\subseteq
\left\{
\sum_t\htr_{t+1}
\left(\delta^*+\frac{\rho_t}{X_t}\right)=+\infty
\right\}
\qquad\text{a.s.}
\]
\item\label{cor-X0:b}
Consider the dynamics~\eqref{eq:recursion:mart} under
Assumptions~\ref{ass:reinforcement} and~\ref{ass:step-size-1b},
and suppose that
\begin{equation}\label{eq:signed-variance-drift}
\Var(\widetilde V_{t+1}\mid\F_t)
=O\!\left(\E(\widetilde V_{t+1}\mid\F_t)\right).
\end{equation}
Then \footnote{\label{footnote}The conditions
are understood starting from a sufficiently large random time (inclusive) after which \(X_t>0\).}
\begin{equation}\label{eq-dicotomy-H}
E_b:=
\liminf_t\{X_t>0\}
\cap
\left\{\limsup_t\frac{\rho_t^-}{X_t}<\delta^*\right\}
\subseteq
\{X^{(\pi)}_t\to0\}
\cup
\{X^{(\pi)}_t\to+\infty\}
\;\text{a.s.}\,,
\end{equation}
and 
\begin{equation}\label{eq-criterion-Hinfty}
E_b\cap \Big\{
\sum_t\frac{\widetilde V_{t+1}^-}{X^{(\pi)}_t}<+\infty
\Big\}\subseteq\{X^{(\pi)}_t\to+\infty\}
\qquad\mbox{a.s.}\,.
\end{equation}
\item\label{cor-X0:c} Consider the dynamics~\eqref{eq:recursion:mart}
under Assumptions~\ref{ass:reinforcement}
and~\ref{ass:step-size-1b}. Let \(E_b\) be the event 
defined in~\eqref{eq-dicotomy-H}. Then~\footref{footnote}
\begin{equation}\label{eq:relative-variance-Hinfty}
E_c:=E_b\cap
\left\{
\sum_t
\frac{\Var(\widetilde V_{t+1}\mid\F_t)}
     {(X_t^{(\pi)})^2}<+\infty
\right\}
\subseteq
\{X_t^{(\pi)}\to+\infty\}
\qquad\text{a.s.}\,.
\end{equation}
\end{enumerate}
\end{theorem}
Note that, for the general case \ref{cor-X0:b}, if $P(E_b)=1$, then, with probability one, only the two alternatives, $\{X^{(\pi)}_t\to 0\}$ or $\{X^{(\pi)}_t\to +\infty\}$, are possible.
In particular, on $\{\sum_t{\widetilde V_{t+1}^-}/{X^{(\pi)}_t}<+\infty\}$, only the latter is possible.
An example showing that $X^{(\pi)}_t\to 0$ may occur with strictly positive probability on $E_b\cap \{\sum_t{\widetilde V_{t+1}^-}/{X^{(\pi)}_t}=+\infty\}$ is presented in Subsection~\ref{sec:H_t_converge_zero}.

\section{Proofs}\label{sec:proofs}
The proofs are ordered by dependence. In the whole sequel, when we work on a given event on which the stated
conditions hold eventually with probability one, all the arguments work after the standard localization to that
event. Moreover, as already said, on the event $\{X^{(\pi)}_t\to+\infty\}$, all the quantities and the sums containing 
$1/X^{(\pi)}_t$ (and $1/X_t$) are understood starting from a sufficiently large random time (inclusive) after which we have $X^{(\pi)}_t>0$ (and so $X_t>0$); 
this changes only finitely many terms. 

\subsection{Proofs of the main results}
We premise the following technical result.
\begin{lemma}[Logarithmic comparison]
\label{lem:multiplicative-log}
Let \((Z_t)_t\) be a strictly positive adapted process
such that
\[
Z_{t+1}=Z_t(1+\mu_{t+1}+\Delta \mathcal M^Z_{t+1}),
\]
where \(\mu_{t+1}\) is \(\F_t\)-measurable and
\((\Delta\mathcal M^Z_{t+1})_t\) are the increments
of a locally square-integrable martingale \footnote{There exists a localizing sequence of stopping times $(\tau_n)_{n \ge 1}$ increasing to $+\infty$ a.s. 
such that, for every $n \ge 1$, the stopped process $({\mathcal M}^{Z})^{|\tau_n}$ is a square-integrable martingale.}.
Suppose that
\[
\sum_{t}\mu_{t+1}^2<+\infty,
\qquad
\sum_{t}
\E( (\Delta \mathcal M^Z_{t+1})^2\mid\F_t)<+\infty
\qquad\text{a.s.}\,.
\]
Then, for a real random variable \(L_Z\),
\[
\ln Z_T-\sum_{t\le T-1}\mu_{t+1}
\stackrel{a.s.}\longrightarrow L_Z.
\]
\end{lemma}
\begin{proof} 
First, notice that the local martingale
\(\sum_{t<T}\Delta\mathcal M^Z_{t+1}\) has finite
predictable quadratic variation by hypothesis. Therefore Localized
Martingale Convergence Theorem \cite[Theorem~2.3.4, p.~46]{lapeyre2024montecarlo}
gives the almost sure convergence of $\mathcal{M}^Z_T$ to a real random variable.

Let us now focus on proving 
\begin{equation}\label{eq:M_Z square convergence}
\sum_t(\Delta\mathcal M^Z_{t+1})^2<+\infty
\qquad\text{a.s.}
\end{equation}
To this end, set
\[
v_{t+1}:=
\E\bigl((\Delta\mathcal M^Z_{t+1})^2\mid\F_t\bigr),
\qquad
Q_T:=
\sum_{t<T}\bigl((\Delta\mathcal M^Z_{t+1})^2\wedge1\bigr).
\]
The increments of \(Q_T\) lie in \([0,1]\); hence
\[
\Var(Q_{t+1}-Q_t\mid\F_t)
\leq \E(Q_{t+1}-Q_t\mid\F_t)
\leq v_{t+1}\quad\mbox{a.s.}\,.
\]
Since \(\sum_t v_{t+1}<+\infty\) a.s.,
Lemma~\ref{lem:compensated-increments}\ref{lem:compensated-increments:i}
shows that \(Q_T\) converges a.s. to a real random variable. Thus, since by the convergence of ${\mathcal M}^Z_T$ only finitely many
squared increments $(\Delta \mathcal{M}^Z_{t+1})^2$ exceed \(1\), we get~\eqref{eq:M_Z square convergence}.
Thus, with
\(U_{t+1}:=\mu_{t+1}+\Delta \mathcal M^Z_{t+1}\), we have
\(\sum_tU_{t+1}^2<+\infty\) a.s. In particular,
\(U_{t+1}\stackrel{a.s.}\to 0\). Since \(1+U_{t+1}>0\), eventually
\(
\bigl|\ln(1+U_{t+1})-U_{t+1}\bigr|
\leq 2U_{t+1}^2.
\)
The logarithmic remainders therefore form a convergent series.
Taking the sum in 
\[
\ln Z_{t+1}-\ln Z_t
=
\mu_{t+1}+\Delta \mathcal M^Z_{t+1}
+\bigl[\ln(1+U_{t+1})-U_{t+1}\bigr]
\]
proves the assertion.
\end{proof}

\begin{proof}[Proof of Theorem~\ref{thm:main-positive-all-maringale-noise} and Theorem~\ref{thm:positive-all_second-maringale-noise}]
By Assumption~\ref{ass:reinforcement},
$
    Y_{t+1}
    =
    \delta^*X_t+\rho_t+\Delta M^Y_{t+1}.
$
Substituting this decomposition into \eqref{eq:recursion:mart} gives
$$
    X_{t+1}
    =
    a_{t+1}X_t
    +
    r_{t+1}\rho_t
    +
    r_{t+1}\Delta M^Y_{t+1}
    +
    r_{t+1}X_t\Delta N_{t+1},
$$
where \(a_{t+1}\) is defined in
\eqref{eq:linear-coefficient:mart2}. We now work on \(\{X^{(\pi)}_t\to+\infty\}\) and consider the recursion at a sufficiently large 
$t$ from which we have $X^{(\pi)}_t>0$ (and so $X_t>0$). 
Since
$q_{t+1}^{(\delta^*)}/q_t^{(\delta^*)}=1/{a_{t+1}}$, 
we obtain
\begin{align}
X_{t+1}^*
&=
X_t^*
\left(
1+
\mu_{t+1}+ 
\Delta\mathcal M_{t+1}
\right),
\label{eq:relative-normalized-recursion:mart2}
\end{align}
where
\begin{equation}\label{eq:predictable-relative-drift:mart2}
\mu_{t+1}
:=
\frac{r_{t+1}\rho_t}{a_{t+1}X_t}
=
\frac{\htr_{t+1}}
{1+\delta^*\htr_{t+1}}
\frac{\rho_t}{X_t},
\end{equation}
and \(\Delta\mathcal M_{t+1}\) is the delta-martingale given by
\eqref{eq:relative-martingale-increment:mart2}.     
Since Assumption~\ref{ass:reinforcement}, and so the consequent bound~\eqref{eq:bound-rho-meno}, 
from~\eqref{eq:predictable-relative-drift:mart2} we get
$$
    {|\mu_{t+1}|}
    \leq
    \htr_{t+1}
    \frac{\rho_t^+}{X_t}
    +
    \delta^*\htr_{t+1}.
$$
It follows that
$$
    {\mu_{t+1}^2}
    \leq
    2\left(
        \htr_{t+1}\frac{\rho_t^+}{X_t}
    \right)^2
    +
    2(\delta^*)^2(\htr_{t+1})^2.
$$
Hence Assumptions~\ref{ass:step-size-2} and
\ref{ass:quadratic-relative-perturbation} imply
\(\sum_t 
{\mu_{t+1}^2}<+\infty
\)
a.s. 
Then Assumption~\ref{ass:relative-noise:mart2} 
allow us to apply
Lemma~\ref{lem:multiplicative-log} with
$Z_t=X_t^*$, $\mu_{t+1}=
{\mu_{t+1}}$ and $\Delta\mathcal M^Z_{t+1}=\Delta\mathcal M_{t+1}$. 
Consequently, for a real random variable \(L'\),
\begin{equation}\label{eq:log-normalized-decomposition:mart2}
\ln X_{T+1}^*
-
\sum_{t\leq T} {\mu_{t+1}} 
\stackrel{a.s.}\longrightarrow L'
\qquad\text{on }\{X^{(\pi)}_t\to+\infty\}.
\end{equation}
It remains to compare the predictable sum in
\eqref{eq:log-normalized-decomposition:mart2} with \(\Gamma_T\). 
By the definition of \(\Gamma_T\) in~\eqref{eq:Gamma} and
{$\mu_{t+1}$} 
in~\eqref{eq:predictable-relative-drift:mart2}, we get
\begin{equation}
\htr_{t+1}\frac{\rho_t}{X_t}
-
{\mu_{t+1}}
=
\frac{
\delta^*(\htr_{t+1})^2
}{
1+\delta^*\htr_{t+1}
}
\frac{\rho_t}{X_t}
=
\frac{
\delta^*(\htr_{t+1})^2
}{
1+\delta^*\htr_{t+1}
}
\frac{\rho_t^+}{X_t} -
\frac{
\delta^*(\htr_{t+1})^2
}{
1+\delta^*\htr_{t+1}
}
\frac{\rho_t^-}{X_t}
.
\label{eq:Gamma-drift-difference:mart2}
\end{equation}
For the positive part, the Cauchy--Schwarz inequality gives
\begin{multline*}
\sum_t
\frac{
\delta^*(\htr_{t+1})^2
}{
1+\delta^*\htr_{t+1}
}
\frac{\rho_t^+}{X_t}
\leq
\delta^*
\sum_t
\htr_{t+1}
\left(
\htr_{t+1}
\frac{\rho_t^+}{X_t}
\right)
\\
\leq
\delta^*
\left(
\sum_t(\htr_{t+1})^2
\right)^{1/2}
\left[
\sum_t
\left(
\htr_{t+1}
\frac{\rho_t^+}{X_t}
\right)^2
\right]^{1/2},
\end{multline*}
which is a.s. finite by Assumptions~\ref{ass:step-size-2} and
\ref{ass:quadratic-relative-perturbation}.
For the negative part, by Assumption~\ref{ass:reinforcement} and so the bound~\eqref{eq:bound-rho-meno}, 
we have 
$$
    \sum_t
    \frac{
        \delta^*(\htr_{t+1})^2
    }{
        1+\delta^*\htr_{t+1}
    }
    \frac{\rho_t^-}{X_t}
    \leq
    (\delta^*)^2
    \sum_t(\htr_{t+1})^2,
$$
which is again a.s. finite by Assumptions~\ref{ass:step-size-2}. 
Therefore we obtain 
\begin{equation}\label{eq:Gamma-drift-limit:mart2}
\sum_{t\leq T}{\mu_{t+1}} 
-\Gamma_T
\stackrel{a.s.}\longrightarrow \widetilde{L}\,,
\end{equation}
for some real random variable \(\widetilde{L}\).
Combining \eqref{eq:log-normalized-decomposition:mart2}
and \eqref{eq:Gamma-drift-limit:mart2}, we obtain
\[
\ln X_{T+1}^*-\Gamma_T
\stackrel{a.s.}\longrightarrow L'+\widetilde L
\qquad\text{on }\{X^{(\pi)}_t\to+\infty\}.
\]
This proves \eqref{eq:logX-Gamma-additive}. Taking lower and
upper limits proves the two stated equivalences and concludes
the proof of Theorem~\ref{thm:main-positive-all-maringale-noise}. 
\\
\indent Assume now Assumption~\ref{ass:lim_sup2}.
Theorem~\ref{thm:finite-limit}, applied to
\eqref{eq:recursion:mart}, gives the a.s. convergence of $X_t^*$ toward a real random variable $X^*_\infty$ and
\begin{equation}\label{eq:rho-absolute}
        \sum_{t}q_{t+1}^{(\delta^*)}r_{t+1}|\rho_t|< +\infty
        \qquad\text{a.s.}
\end{equation}
If $X^*_\infty>0$ a.s.\,, then 
$1/X_t^*=O(1)$ and
\(
q_t^{(\delta^*)}\htr_{t+1}
=q_{t+1}^{(\delta^*)}r_{t+1}
\bigl(1+\delta^*\htr_{t+1}\bigr),
\)
so \eqref{eq:rho-absolute} proves
\(
\sum_t\htr_{t+1}|\rho_t|/X_t<+\infty
\) a.s.\,.
Thus \ref{pos-I2} implies \ref{pos-VI2}.  The absolute
 summability~\ref{pos-VI2} makes $\Gamma_T$ a.s. convergent and implies
Assumption~\ref{ass:quadratic-relative-perturbation}, so
\ref{pos-VI2} implies \ref{pos-V2}.  Finally,
\ref{pos-V2} permits the logarithmic comparison proved above in Theorem~\ref{thm:main-positive-all-maringale-noise},
which gives $\liminf_TX_T^*>0$ and hence
\ref{pos-I2}.  Therefore
\[
\ref{pos-I2}\Longleftrightarrow
\ref{pos-V2}\Longleftrightarrow
\ref{pos-VI2} 
\]
and also Theorem~\ref{thm:positive-all_second-maringale-noise} is proven.
\end{proof}

For the proof of Theorem~\ref{rem:conditional-variance-implies-relative-noise:mart2}, we use the following result. 

\begin{theorem}[Logarithmic asymptotics]
\label{thm:log-compensator_comp-relative}
Let  $V_{t+1}\geq 0$ and suppose that
Assumption~\ref{ass:moment} holds. Set $H_{t+1}=H_t+V_{t+1}$ and assume  
\begin{equation}\label{eq:log_inc_mean_o}
\frac{\E(V_{t+1}\mid\F_t)}{H_t}\stackrel{a.s.}\longrightarrow 0
\qquad\text{on }\{H_t\to+\infty\}.
\end{equation}
Then \footnote{All the quantities 
are understood starting from a sufficiently large random time (inclusive) after which \(H_t\geq 1\).}
\begin{equation}\label{eq:log-compensator}
\ln H_{T+1}    \stackrel{a.s.} \sim  \sum_{t\leq T}\frac{\E(V_{t+1}\mid\F_t)}{H_t}
   \qquad\text{on } \{H_t\to+\infty\},
\end{equation}
and
\begin{equation}\label{eq:mu-H2-summable}
    \sum_t\frac{\E(V_{t+1}\mid\F_t)}{H_t^2}<+\infty\qquad\text{a.s. on } \{H_t\to+\infty\}.
\end{equation}
Consequently,
\begin{equation}\label{eq:relative-variance-summable-basic}
    \sum_t
    \frac{\Var(V_{t+1}\mid\F_t)}{H_t^2}
    <+\infty
    \qquad\text{a.s. on }\{H_t\to+\infty\}.
\end{equation}
\end{theorem}
\begin{proof}
Work on \(\{H_t\to+\infty\}\), starting from a finite stopping time (inclusive) after which
\(H_t\geq 1\). Write
\[
W_t:=\frac{\E(V_{t+1}\mid\F_t)}{H_t},
\qquad
L_t:=
\E\left(
 \ln\frac{H_{t+1}}{H_t}\big|\F_t
\right)=\E\left(
 \ln\bigl(1+\frac{V_{t+1}}{H_t}\bigr)\mid\F_t
\right).
\]
Assumption~\ref{ass:moment} and the inequalities
\(x-x^2/2\leq\ln(1+x)\leq x\), for \(x\geq0\),
give
\[
0\leq W_t- L_t
\leq
\frac12\left(
 \frac{O(\E(V_{t+1}\mid\F_t))}{H_t^2}+W_t^2
\right)
=
\frac{W_t}{2}
\left(\frac{O(1)}{H_t}+W_t\right).
\]
Since  we are working on \(\{H_t\to+\infty\}\) and  we have \(W_t\stackrel{a.s.}\to 0\)
by \eqref{eq:log_inc_mean_o}, from the above inequalities we obtain \(L_t=W_t(1+o(1))\) on \(\{W_t>0\}\) and \(L_t=W_t=0\) otherwise. 
\\
\indent Now, 
we apply Lemma~\ref{lem:compensated-increments} to the
non-negative (starting from a sufficiently large time) adapted process  $Z_t:=\ln(H_t)$. 
 Indeed, since the conditional expectation minimizes the mean squared error among all \(\F_t\)-measurable random variables and, on
\([0,+\infty)\), the function \(x\mapsto\ln(1+x)\) is Lipschitz with constant equal to $1$, we have on $\{H_t\to +\infty\}$
 \begin{multline*}
  \Var\Bigl(
        \ln\Bigl(1+\frac{V_{t+1}}{H_t}\Bigr)
        \big|\F_t\Bigr)\leq
\E\biggl( \Bigl( \ln\Bigl(1+\frac{V_{t+1}}{H_t}\Bigr)-\ln\Bigl(1+W_t\Bigr) \Bigr)^2 \Big| {\F}_t \biggr)
  \\
  \leq
        \E\biggl( \Bigl( \frac{V_{t+1}}{H_t}-W_t\Bigr)^2\Big| \F_t\biggr)
  = \Var\Bigl(\frac{V_{t+1}}{H_t}\big|\F_t\Bigr)=\frac{\Var(V_{t+1}\mid\F_t)}{H_t^2} \,,
 \end{multline*}
 that is, by Assumption~\ref{ass:moment} and the fact that $W_t=O(L_t)$, we have on $\{H_t\to +\infty\}$
\begin{equation*}
\begin{split}
\Var(\Delta Z_{t+1}\mid\F_t)
\leq
\frac{\Var(V_{t+1}\mid\F_t)}{H_t^2}
&=O\Bigl(\frac{\E(V_{t+1}\mid\F_t)}{H_t^2}\Bigr)
=
\frac{O(W_t)}{H_t}\\
&= O(L_t)= O(\E(\Delta Z_{t+1}\mid\F_t)),
\end{split}
\end{equation*}
as required by the above mentioned lemma. Therefore,
Let \(A_T:=\sum_{t\leq T-1} L_t\).
If \(A_\infty<+\infty\), part \ref{lem:compensated-increments:i} of
Lemma~\ref{lem:compensated-increments} would imply that
\(Z_T=\ln H_T\) has a finite limit, contrary to the fact that we are working on 
\(\{H_t\to+\infty\}\). Thus \(A_T\stackrel{a.s.}\to+\infty\), and
part \ref{lem:compensated-increments:ii} gives 
$\ln H_T\stackrel{a.s.}\sim A_T$. Moreover, as stated above, \(L_t=W_t(1+o(1))\) on \(\{W_t>0\}\) and \(L_t=W_t=0\) otherwise. Since \(A_T=\sum_{t\leq T-1} L_t\stackrel{a.s.}\to +\infty\), the weighted Ces\`aro argument gives
\[
\ln H_T \stackrel{a.s.}\sim 
A_T = 
\sum_{t\leq T-1} L_t
\stackrel{a.s.}\sim\sum_{t\leq T-1} W_t = \sum_{t\leq T-1} \frac{\E(V_{t+1}\mid\F_t)}{H_t}\,,
\]
that proves
\eqref{eq:log-compensator}. \\ \indent Finally, apply part \ref{lem:compensated-increments:iii} of the same~Lemma~\ref{lem:compensated-increments}
to the original process \(Z_t=H_t\). Since 
\(\E(V_{t+1}\mid\F_t)/H_t=W_t\stackrel{a.s.}\to 0\), we get on \(\{H_t\to +\infty\}\)  
\[
\sum_t\frac{\E(V_{t+1}\mid\F_t)}{H_t^2}<+\infty
\qquad\mbox{and}\qquad 
\sum_t
\frac{\Var(V_{t+1}\mid\F_t)}{H_t^2}
<+\infty\quad\mbox{a.s.},
\]
which are \eqref{eq:mu-H2-summable} and
\eqref{eq:relative-variance-summable-basic}.
\end{proof}

\begin{proof}[Proof of Theorem~\ref{thm:positive-all_second}] 
Work on $\{X^{(\pi)}_t\to +\infty\}$. 
For the basic recursion, by the definition of $\Delta\mathcal M_{t+1}$ in ~\eqref{eq:relative-martingale-increment:mart2} we have 
\[
\Delta\mathcal M_{t+1}
=
\frac{V_{t+1}-\E(V_{t+1}\mid\F_t)}
     {X^{(\pi)}_t(1+\delta^*\htr_{t+1})}.
\]
Moreover, by Assumption~\ref{ass:reinforcement}, 
we can write 
\[
0\leq
\frac{\E(V_{t+1}\mid\F_t)}{X^{(\pi)}_t}
=
\htr_{t+1}
\Bigl(\delta^*+\frac{\rho_t}{X_t}\Bigr)
\leq
\delta^*\htr_{t+1}
+\htr_{t+1}\frac{\rho_t^+}{X_t}.
\]
Hence, Assumptions~\ref{ass:step-size-2} and
\ref{ass:quadratic-relative-perturbation} imply
\(
\sum_t
(
{\E(V_{t+1}\mid\F_t)}/{X^{(\pi)}_t}
)^2
<+\infty.
\)
In particular,
\(
{\E(V_{t+1}\mid\F_t)}/{X^{(\pi)}_t}\to 0,
\)
so condition~\eqref{eq:log_inc_mean_o} in
Theorem~\ref{thm:log-compensator_comp-relative} holds. Therefore, by using~\eqref{eq:log_inc_mean_o} and Assumption~\ref{ass:moment} we can apply Theorem~\ref{thm:log-compensator_comp-relative} to get 
\eqref{eq:relative-variance-summable-basic} that gives
\[
\begin{aligned}
\sum_t
\E\bigl((\Delta\mathcal M_{t+1})^2\mid\F_t\bigr)
=
\sum_t
\frac{\Var(V_{t+1}\mid\F_t)}
     {\bigl(X^{(\pi)}_t\bigr)^2(1+\delta^*\htr_{t+1})^2}  
\leq
\sum_t\frac{\Var(V_{t+1}\mid\F_t)}{\bigl(X^{(\pi)}_t\bigr)^2}
<+\infty.
\end{aligned}
\]
Thus Assumption~\ref{ass:relative-noise:mart2} holds. The conclusions
now follow from
Theorems~\ref{thm:main-positive-all-maringale-noise}
and~\ref{thm:positive-all_second-maringale-noise}.
\end{proof}

\subsection{Proofs of the results useful for applications}

\begin{proof}[Proof of Theorem~\ref{prop:absolute-relative-criterion}]
On \(G\), consider the recursion~\eqref{eq:relative-normalized-recursion:mart2}, that is 
\[
 X_{t+1}^*
 =
 X_t^*\bigl(1+
 {\mu_{t+1}}+
 \Delta\mathcal M_{t+1}\bigr),
\]
starting from  a sufficiently large random time (inclusive) after which we have $X^{(\pi)}_t>0$ and so  \(X_t>0\).
The expression~\eqref{eq:predictable-relative-drift:mart2} for 
\(\mu_{t+1}\) and the 
absolute relative-perturbation condition on \(G\) give
\(
 \sum_{t} 
 {|\mu_{t+1}|}<+\infty 
\)
and so 
$ \sum_{t} 
{\mu_{t+1}^2} <+\infty$ on $G$. 
Moreover, Assumption~\ref{ass:relative-noise:mart2} gives
\[
 \sum_{t}
 \E\bigl(
  (\Delta\mathcal M_{t+1})^2\mid\F_t
 \bigr)<+\infty
 \qquad\text{a.s. on }G.
\]
Lemma~\ref{lem:multiplicative-log} therefore yields
\[
 \ln X_T^*
 -
 \sum_{t\leq T-1} {\mu_{t+1}} 
 \stackrel{a.s.}\longrightarrow L^*
 \qquad\text{on } G,
\]
for a real random variable \(L^*\). Hence, since $\sum_t |\mu_{t+1}|<+\infty$ on $G$, we can affirm that on $G$ the process  \(\ln X_T^*\)
converges a.s. to a real random variable, that is on $G$ the process $X^*_T$ 
converges a.s. to a strictly positive real random variable.
\end{proof}

\begin{proof}[Proof of Theorem~\ref{cor-X0}]
\ref{cor-X0:a}
Work on \(\bigcup_t\{X_t>0\}\). On this event, the stopping time 
\(
\tau:=\inf\{t\geq 0:X_t>0\} 
\)
takes values in $\mathbb N$  and
\(X^{(\pi)}_t\geq X^{(\pi)}_\tau>0\) for every \(t\geq\tau\). 
Assumption~\ref{ass:reinforcement} gives, for \(t\geq \tau\),
\[
\E(V_{t+1}\mid\F_t)
=\htr_{t+1}X^{(\pi)}_t
  \left(\delta^*+\frac{\rho_t}{X_t}\right)
\geq X^{(\pi)}_{\tau}\htr_{t+1}
  \left(\delta^*+\frac{\rho_t}{X_t}\right).
\]
The quantities in parentheses are non-negative because of~\eqref{eq:bound-rho-meno}. 
Thus, on the event $E_a$,
we have 
\(
\sum_{t}\E(V_{t+1}\mid\F_t)=+\infty.
\)
Assumption~\ref{ass:moment} allows us to apply
Lemma~\ref{lem:compensated-increments}\ref{lem:compensated-increments:ii} to
\(X^{(\pi)}_{t+1}-X^{(\pi)}_t=V_{t+1}\). Consequently, on $E_a$
\[
X^{(\pi)}_{T+1}
\stackrel{a.s.}\sim
\sum_{t\leq T}\E(V_{t+1}\mid\F_t)
\stackrel{a.s.}\longrightarrow+\infty,
\]
which proves \eqref{eq-criterio-Hinfty2}. 
For the final inclusion in \ref{cor-X0:a}, use
\[
\delta^*+\frac{\rho_t}{X_t}
=
\frac{\rho_t^+}{X_t}
+\left(\delta^*-\frac{\rho_t^-}{X_t}\right).
\]
On either event in that inclusion, the expression on the
left is eventually bounded below by a strictly positive
random quantity. Assumption~\ref{ass:step-size-1b} then makes
its product with \(\htr_{t+1}\) have a divergent sum.\\
\ref{cor-X0:b}
Define 
\(
A_T:=\sum_{t=0}^{T-1}
       \E(\widetilde V_{t+1}\mid\F_t).
\)
We apply Lemma~\ref{lem:compensated-increments} to the adapted
process \(X^{(\pi)}_t\) and its compensator $A_T$. Its conditional
drift is non-negative, since
\(\E(\widetilde V_{t+1}\mid\F_t)
 =\E(V_{t+1}\mid\F_t)\geq0\). 
Now work on \(E_b\), and choose an a.s. finite random
time \(\tau\) such that \(X_t>0\) for every \(t\geq\tau\) and
\(\delta^*-\rho_t^-/X_t\geq\eps\) for every \(t\geq\tau\),
for some random \(\eps>0\). Then, for \(t\geq\tau\), 
\begin{equation}\label{eq:signed-effective-drift-proof}
\E(\widetilde V_{t+1}\mid\F_t)
=
\htr_{t+1}X^{(\pi)}_t
\left(\delta^*+\frac{\rho_t}{X_t}\right)
\geq\eps\htr_{t+1}X^{(\pi)}_t
\geq0.
\end{equation}
Apply Lemma~\ref{lem:compensated-increments} to the process
\(X^{(\pi)}_t\), whose increments are \(\widetilde V_{t+1}\);
condition~\eqref{eq:signed-variance-drift} supplies its
conditional-variance bound. 
On \(\{A_T\to+\infty\}\), part \ref{lem:compensated-increments:ii} of the lemma gives
\(X^{(\pi)}_T\stackrel{a.s.}\sim A_T\to+\infty\).
On  \(\{A_\infty<+\infty\}\), part
\ref{lem:compensated-increments:i} gives a finite limit
\(X^{(\pi)}_\infty\). This limit is zero on \(E_b\).
Indeed, if \(X^{(\pi)}_\infty>0\), then
\eqref{eq:signed-effective-drift-proof} implies, eventually,
\[
\E(\widetilde V_{t+1}\mid\F_t)
\geq \frac{\eps X^{(\pi)}_\infty}{2}
       \htr_{t+1}.
\]
Then, Assumption~\ref{ass:step-size-1b} would give
\(A_\infty=+\infty\), a contradiction. This proves
\eqref{eq-dicotomy-H}.
\\
\indent Now work on $E_b \cap \{\sum_t{\widetilde V_{t+1}^-}/{X^{(\pi)}_t}<+\infty\}$.
Enlarge \(\tau\), if needed, so that \(0\leq U_{t+1}<1\) for every \(t\geq\tau\).
Set
\(
U_{t+1}:={\widetilde V_{t+1}^-}/{X^{(\pi)}_t}
\) so that \(0\leq U_{t+1}<1\) eventually and
\(\sum_tU_{t+1}<+\infty\).
Since
\(
X^{(\pi)}_{t+1}
\geq X^{(\pi)}_t(1-U_{t+1}),
\)
iteration starting from  
\(\tau\) yields
\[
X^{(\pi)}_T
\geq X^{(\pi)}_{\tau}
\prod_{t=\tau}^{T-1}(1-U_{t+1})\,,
\]
where $X^{(\pi)}_{\tau}>0$ and  the product has a strictly positive limit. Thus
\(\liminf_T X^{(\pi)}_T>0\), and so, by \eqref{eq-dicotomy-H}, we get~\eqref{eq-criterion-Hinfty}.
\\
\ref{cor-X0:c}
Define 
\[U_t:=\E(\widetilde V_{t+1}\mid\F_t)=\pi_{t+1}r_{t+1}
 \left(\delta^*X_t+\rho_t\right)\geq 0.
\]
Moreover, set \(D_{t+1}:=0\) and $W_t:=0$ when \(X_t^{(\pi)}=0\), and,
when \(X_t^{(\pi)}>0\), set 
\[
D_{t+1}:=
\frac{\widetilde V_{t+1}-U_t}
     {X_t^{(\pi)}+U_t},
\qquad
W_t:=
\E(D_{t+1}^2\mid\F_t)
=
\frac{\Var(\widetilde V_{t+1}\mid\F_t)}
     {(X_t^{(\pi)}+U_t)^2}.
\]
In particular, \(W_t\) is
\(\F_t\)-measurable. For each integer \(n\geq1\), set
\(
D_{t+1}^{[n]}
:=\mathbf1_{\{W_t\leq n\}}D_{t+1}.
\) 
Its cutoff is \(\F_t\)-measurable, and
\[
\E(D_{t+1}^{[n]}\mid\F_t)=0,
\qquad
\E\bigl((D_{t+1}^{[n]})^2\mid\F_t\bigr)
\leq n.
\]
Hence \(\mathcal D_T^{[n]}:=\sum_{t\leq T-1}D_{t+1}^{[n]}\)
is a square-integrable martingale.
\\
\indent Now work on \(E_c\). Choose an a.s. finite random time
\(\tau\) such that \(X_t>0\), and  so \(X_t^{(\pi)}>0\), for every \(t\geq\tau\).
Thus, 
\begin{equation}\label{eq:bracket-D}
\sum_{t\geq\tau}W_t
\leq
\sum_{t\geq\tau}
\frac{\Var(\widetilde V_{t+1}\mid\F_t)}
     {(X_t^{(\pi)})^2}
<+\infty\qquad\mbox{a.s. on  } E_c\,.
\end{equation}
Consequently, for every fixed \(n\), the bracket of the martingale 
\(\mathcal D^{[n]}\) is a.s. finite on \(E_c\): the terms before
\(\tau\) are finite in number and bounded by \(n\).
The Localized Martingale Convergence Theorem \cite[Theorem~2.3.4, p.~46]{lapeyre2024montecarlo} therefore gives the a.s. 
convergence of \(\mathcal D_T^{[n]}\) on \(E_c\). To prove 
\begin{equation}\label{eq:square-sum}
\sum_t(D_{t+1}^{[n]})^2<+\infty
\qquad\text{a.s. on } E_c\,,
\end{equation}
set 
\[
Q_T^{[n]}
:=
\sum_{t\leq T-1}\bigl((D_{t+1}^{[n]})^2\wedge1\bigr).
\]
Its increments lie in \([0,1]\), and their conditional expectations 
are bounded by $\E\bigl((D_{t+1}^{[n]})^2\mid\F_t\bigr)$, whose series is a.s.\ finite on $E_c$ by \eqref{eq:bracket-D}. 
Lemma~\ref{lem:compensated-increments}\ref{lem:compensated-increments:i}
thus shows that \(Q_T^{[n]}\) converges a.s. on \(E_c\).
Since 
\(D_{t+1}^{[n]}\stackrel{a.s.}\to 0\), we obtain~\eqref{eq:square-sum}. 
Because \(\sum_{t\geq\tau} W_t<+\infty\) a.s. on $E_c$, we have 
$$
E_c\subseteq \bigcup_n \{D_{t+1}^{[n]}=D_{t+1}\;\forall t\geq \tau\}\qquad\mbox{a.s.}\,.
$$
As a consequence, the a.s. convergences for the series involving $D_{t+1}^{[n]}$ yield 
\[
\sum_{t\geq\tau} D_{t+1}\ \text{a.s. converges on } E_c,
\qquad
\sum_{t\geq\tau}D_{t+1}^2<+\infty
\quad\text{a.s. on }E_c.
\]
For \(t\geq\tau\), set
\(
\widetilde{U}_t:={U_t}/{X_t^{(\pi)}}\). 
Then
\[
X_{t+1}^{(\pi)}
=X_t^{(\pi)}(1+\widetilde{U}_t)(1+D_{t+1}).
\]
Since $X_t^{(\pi)}>0$ for $t\geq \tau$, we have \(1+D_{t+1}>0\) for $t\geq \tau$.
The series convergences just proved and Taylor's formula show that
\(\sum_{t\geq\tau}\ln(1+D_{t+1})\) a.s. converges. 
\\
\indent Finally, by the definition of \(E_b\), we may enlarge \(\tau\) so that,
for some random \(\eps>0\),
\(\widetilde{U}_t\geq\eps \htr_{t+1}\) whenever \(t\geq\tau\).
Since, for every \(x\geq 0\), 
\(
\ln(1+\eps x)
\geq\ln(1+\eps)\min\{x,1\}
\), 
Assumption~\ref{ass:step-size-1b} implies
\[
\sum_{t\geq\tau}\ln(1+\widetilde{U}_t)
\geq
\sum_{t\geq\tau}\ln(1+\eps \htr_{t+1})
=+\infty.
\]
Taking logarithms of the above recursion for \(X_t^{(\pi)}\)
now gives \(\ln X_T^{(\pi)}\stackrel{a.s.}\to+\infty\) on \(E_c\), proving~\eqref{eq:relative-variance-Hinfty}.
\end{proof}

\section{Complements on the basic stochastic recursion}\label{sec:lower}

The exact-normalization theorems require second-order summability.  This
section consider the basic stochastic recursion~\eqref{eq:recursion} and records the weaker first-order identities that remain available
without it and then derives auxiliary criteria used to verify divergence
and perturbation summability.

\begin{proposition}
\label{thm:log-compensator}
Let  $V_{t+1}\geq 0$ and suppose that
Assumption~\ref{ass:moment} holds. Moreover, set $H_{t+1}=H_t+V_{t+1}$.  
Then we have \footnote{All sums involving a denominator \(P_t\) or \(H_t\) are
understood from a sufficiently large time at which that 
denominator is strictly positive.}:
\begin{enumerate}[label=\textup{(\alph*)}]
\item\label{comp-replacement}[Predictable replacement] 
For every adapted non-negative process $(P_t)$, on
$\{\liminf_t P_t>0\}$, we have 
\begin{equation}\label{eq:replacement-convergence}
    \sum_t\frac{V_{t+1}}{P_t}<+\infty\quad\text{a.s.}
    \quad\Longleftrightarrow\quad
    \sum_t\frac{\E(V_{t+1}\mid\F_t)}{P_t}<+\infty \quad\text{a.s.}
\end{equation}
and, if these series diverge,
\begin{equation}\label{eq:replacement-ratio}
    \sum_{t\leq T}\frac{V_{t+1}}{P_t}
    \stackrel{a.s.}\sim
    \sum_{t\leq T}\frac{\E(V_{t+1}\mid\F_t)}{ P_t}.
\end{equation}
\item\label{comp-additive}[Additive logarithmic comparison]
On
\begin{equation}\label{eq:relative-mu-square}
    { \{H_t\to+\infty\} \cap \Bigl\{}\sum_t
    \left(
        \frac{\E(V_{t+1}\mid\F_t)}{H_t}
    \right)^2
    <+\infty {\Bigr\}}
\end{equation}
we have
\begin{equation}\label{eq:additive-log-comparison}
    \ln H_{T+1}
    -
    \sum_{t\leq T}
    \frac{\E(V_{t+1}\mid\F_t)}{H_t}
    \stackrel{a.s.}\longrightarrow L_H
\end{equation}
for a real random variable \(L_H\).
\end{enumerate}
\end{proposition}

\begin{proof}
For each \(n\geq1\), set
\(P_t^{[n]}:=P_t\vee n^{-1}\), and define
\[
S_T^{[n]}:=\sum_{t<T}\frac{V_{t+1}}{P_t^{[n]}},
\qquad
A_T^{[n]}:=\sum_{t<T}
 \frac{\E(V_{t+1}\mid\F_t)}{P_t^{[n]}}.
\]
Assumption~\ref{ass:moment} gives for each fixed $n$
\[
\Var\Bigl(
 \frac{V_{t+1}}{P_t^{[n]}}\big|\F_t\Bigr)
=
\frac{\Var(V_{t+1}\mid\F_t)}{(P_t^{[n]})^2}
\leq
n\frac{\Var(V_{t+1}\mid\F_t)}{P_t^{[n]}}
=
O\biggl(
 \frac{\E(V_{t+1}\mid\F_t)}{P_t^{[n]}}
\biggr).
\]
Apply Lemma~\ref{lem:compensated-increments} to
\(S^{[n]}\), for each fixed \(n\). On
\(\{\liminf_tP_t>0\}\), some \(n\) satisfies
\(P_t^{[n]}=P_t\) for all sufficiently large time.
The original and truncated sums therefore differ by only
finitely many terms, giving both conclusions.

\bigskip

For \ref{comp-additive}, 
work on \eqref{eq:relative-mu-square}
starting from a finite stopping time (inclusive) after which
\(H_t\geq 1\).
The second event of
\eqref{eq:relative-mu-square} implies
\(\E(V_{t+1}\mid\F_t)=o(H_t)\). Hence
Lemma~\ref{lem:compensated-increments}\ref{lem:compensated-increments:iii},
applied to \(H_{t+1}=H_t+V_{t+1}\), gives on \eqref{eq:relative-mu-square}
\[
\sum_{t}
\frac{\Var(V_{t+1}\mid\F_t)}{H_t^2}
<+\infty
\qquad\text{a.s.}.
\]
Now write
\[
H_{t+1}
=
H_t\left(
 1+\frac{\E(V_{t+1}\mid\F_t)}{H_t}
 +\frac{V_{t+1}-\E(V_{t+1}\mid\F_t)}{H_t}
\right).
\]
The predictable relative increments are square-summable
by \eqref{eq:relative-mu-square}. The centered relative
increments form a martingale-difference sequence and
satisfy
\[
\sum_t
\E\biggl(
 \Bigl(\frac{V_{t+1}-\E(V_{t+1}\mid\F_t)}{H_t}\Bigr)^2
 \Big|\F_t
\biggr)
=
\sum_t
\frac{\Var(V_{t+1}\mid\F_t)}{H_t^2}
<+\infty.
\]
Also \(H_{t+1}/H_t\geq1\). Thus
Lemma~\ref{lem:multiplicative-log}, applied to \(H_t\)
with \(\mu_{t+1}=\E(V_{t+1}\mid\F_t)/H_t\),
yields
\[
\ln H_{T+1}
-
\sum_{t\leq T}\frac{\E(V_{t+1}\mid\F_t)}{H_t}
\stackrel{a.s.}\longrightarrow L_H
\]
for a suitable real random variable \(L_H\). This is
\eqref{eq:additive-log-comparison}.
\end{proof}

 The following deterministic product expansions and the exact
conditional-mean identity \eqref{eq:EV-deltaR} below, give a first-order
relation between $\ln X_T^*$ and $\Gamma_T$ without any second-order summability 
assumption or boundedness condition on $(X^*_t)$. Specifically,  Lemma~\ref{lem:logZ-identity} shows that the normalized liminf and limsup of these two quantities coincide. 
Theorem~\ref{thm:subexp}, together with the predictable replacement in Proposition~\ref{thm:log-compensator}, then rewrites these identities as separate lower and upper criteria. When both sets of criteria hold, the first-order scale is exact: $\Gamma_T=o(s_T)$, $\ln X_T^*=o(s_T)$, and $\ln X^{(\pi)}_T\stackrel{a.s.}\sim\delta^*s_T$. 

\subsubsection*{Deterministic product estimates.}
Under Assumptions~\ref{ass:step-size-1a} and
\ref{ass:step-size-1b}, for every $\delta>0$,
\begin{equation}\label{eq:PQ-first-order}
        \frac{\ln \pi_{T+1}-\ln q_{T+1}^{(\delta)}}{\delta s_{T+1}}
        \longrightarrow 1
\end{equation}
and,  for every $\eps>0$, we have 
\begin{equation}\label{eq:exp-sum}
        \sum_{t}\htr_{t+1}e^{-\eps s_{t+1}}<+\infty .
\end{equation}
Under the additional Assumption~\ref{ass:step-size-2},
\begin{equation}\label{eq:PQ-second-order}
 \delta s_{T+1}
 -\bigl(\ln \pi_{T+1}-\ln q_{T+1}^{(\delta)}\bigr)
 \longrightarrow d_\delta\in[0,+\infty).
\end{equation}
Indeed,
\begin{equation}\label{eq:PQ-log}
        \ln \pi_{T+1}-\ln q_{T+1}^{(\delta)}=
        \sum_{t\leq T}\ln\left(1+\frac{\delta r_{t+1}}{1-p_{t+1}}\right) =
        \sum_{t\leq T}\ln\left(1+\delta \htr_{t+1}\right) .
\end{equation}
Since $\htr_t\longrightarrow 0$, we get \eqref{eq:PQ-first-order}. 
Moreover, for every $\eps>0$, we have 
\[
        e^{-\eps s_t}-e^{-\eps s_{t+1}}
        =e^{-\eps s_{t+1}}(e^{\eps \htr_{t+1}}-1)
        \ge c_\eps \htr_{t+1}e^{-\eps s_{t+1}}
\]
for all sufficiently large $t$ and, consequently, since $s_T\to +\infty$, we obtain \eqref{eq:exp-sum}.  
For the last assertion, use
\(
0\leq \delta x-\ln(1+\delta x)\leq c_\delta x^2
\) 
for all sufficiently small $x\geq0$, and sum with
$x=\htr_{t+1}$.

The connection with the perturbation is supplied by the conditional mean. Indeed, Assumption~\ref{ass:reinforcement} gives the exact decomposition 
\begin{equation}\label{eq:EV-deltaR} 
\sum_{t\leq T}\frac{\E(V_{t+1}\mid\F_t)}{X^{(\pi)}_t} = \sum_{t\leq T} \frac{\pi_{t+1}r_{t+1}(\delta^*X_t+\rho_t)}{X^{(\pi)}_t}  
= \delta^*(s_{T+1}-s_\tau)+\Gamma_T.
\end{equation}
The next lemma combines this identity with the relative logarithmic compensator theorem and the first-order product expansion. Unlike earlier boundedness-based formulations, it controls both the lower and upper logarithmic behavior without assuming that $X_t^*$ has a finite limsup.

\begin{lemma}[Logarithmic key identity]\label{lem:logZ-identity}
Consider the basic recursion~\eqref{eq:recursion}. Let Assumptions~\ref{ass:reinforcement},  \ref{ass:step-size-1a}, \ref{ass:step-size-1b} and \ref{ass:moment} hold,
together with \eqref{eq:log_inc_mean_o}. Then
\begin{equation}\label{eq:key-identity}
    \ln X_{T+1}^*
    =
    \Gamma_T
    +
    o(s_{T+1}+\Gamma_T^+)\qquad\text{on } \{X^{(\pi)}_t\to+\infty\}.
\end{equation}
Consequently, on the event $\{X^{(\pi)}_t\to+\infty\}$, we have 
\begin{equation}\label{eq:key-liminf-limsup}
\begin{split}
    \liminf_T
    \frac{\ln X_{T+1}^*}{s_{T+1}}
    &=
    \liminf_T
    \frac{\Gamma_T}{s_{T+1}}\qquad\text{a.s.},\\
    \limsup_T
    \frac{\ln X_{T+1}^*}{s_{T+1}}
    &=
    \limsup_T
    \frac{\Gamma_T}{s_{T+1}}\qquad\text{a.s.}.
\end{split}
\end{equation}
Moreover, on the event $\{X^{(\pi)}_t\to+\infty\}$, we have 
\begin{equation}\label{eq:H-liminf-limsup}
\begin{split}
    \liminf_T
    \frac{\ln X^{(\pi)}_{T+1}}{s_{T+1}}
    &=
    \delta^*
    +
    \liminf_T
    \frac{\Gamma_T}{s_{T+1}}\qquad\text{a.s.},\\
    \limsup_T
    \frac{\ln X^{(\pi)}_{T+1}}{s_{T+1}}
    &=
    \delta^*
    +
    \limsup_T
    \frac{\Gamma_T}{s_{T+1}}\qquad\text{a.s.}.
\end{split}
\end{equation}
In particular, we have 
\[
    \liminf_T\frac{\Gamma_T}{s_{T+1}}\geq-\delta^*\qquad\text{a.s. on } \{X^{(\pi)}_t\to+\infty\}. 
\]
\end{lemma}

\begin{proof}
We work on $\{X^{(\pi)}_t\to+\infty\}$. Set
\[
    A_T
    :=
    \sum_{t\leq T}
    \frac{\E(V_{t+1}\mid\F_t)}{X^{(\pi)}_t}.
\]
By Assumption~\ref{ass:reinforcement},
\begin{equation}\label{eq:key-proof-A}
A_T
    =
    \delta^*(s_{T+1}-s_\tau)+\Gamma_T.
\end{equation}
Theorem~\ref{thm:log-compensator_comp-relative} gives
$A_T\stackrel{a.s.}\to+\infty$ 
and
\begin{equation}\label{eq:key-proof-logH}
    \ln X^{(\pi)}_{T+1}
    =
    A_T+o(A_T).
\end{equation}
On the other hand, \eqref{eq:PQ-first-order} yields
\begin{equation}\label{eq:key-proof-product}
    \ln\pi_{T+1}
    -
    \ln q_{T+1}^{(\delta^*)}
    =
    \delta^*s_{T+1}+o(s_{T+1}).
\end{equation}
Therefore, using
\(
    \ln X_{T+1}^*
    =
    \ln X^{(\pi)}_{T+1}
    -
    \ln\pi_{T+1}
    +
    \ln q_{T+1}^{(\delta^*)},
\)
relations \eqref{eq:key-proof-A}--\eqref{eq:key-proof-product} give
\begin{align*}
    \ln X_{T+1}^*
    &=
    A_T+o(A_T)
    -
    \delta^*s_{T+1}
    +
    o(s_{T+1})\\
    &=
    \Gamma_T
    -\delta^*s_\tau
    +
    o(A_T)
    +
    o(s_{T+1}).
\end{align*}
Now $A_T\geq0$, and from \eqref{eq:key-proof-A},
$\Gamma_T\geq-\delta^*s_{T+1}$. 
Moreover,
\[
    0\leq A_T
    =
    \delta^*s_{T+1}-\delta^*s_\tau+\Gamma_T
    \leq
    \delta^*s_{T+1}+\Gamma_T^+.
\]
Hence
\(
    o(A_T)+o(s_{T+1})
    =
    o(s_{T+1}+\Gamma_T^+),
\) 
which proves \eqref{eq:key-identity}.

To obtain \eqref{eq:key-liminf-limsup}, note that $s_{T+1}\to+\infty$, and \eqref{eq:key-identity} can be written as
\begin{equation}\label{eq:key-scaled}
    \frac{\ln X_{T+1}^*}{s_{T+1}}
    =
    \frac{\Gamma_T}{s_{T+1}}+o\Bigl(1+\frac{\Gamma_T^+}{s_{T+1}}\Bigr).
\end{equation}
Notice also that ${\Gamma_T}/{s_{T+1}}\geq-\delta^*$.

Consider any subsequence. If $({\Gamma_T}/{s_{T+1}})$ is bounded above along that
subsequence, then \eqref{eq:key-scaled} gives
\[
    \frac{\ln X_{T+1}^*}{s_{T+1}}-\frac{\Gamma_T}{s_{T+1}}\stackrel{a.s.}\longrightarrow 0.
\]
If instead ${\Gamma_T}/{s_{T+1}}\to+\infty$ along a further subsequence, then
\[
    \frac{{\ln X_{T+1}^*}/{s_{T+1}}}{{\Gamma_T}/{s_{T+1}}}\stackrel{a.s.}\longrightarrow 1,
\]
and therefore ${\ln X_{T+1}^*}/{s_{T+1}}\stackrel{a.s.}\to+\infty$ along the same subsequence. Thus, the two sequences 
$({\ln X_{T+1}^*}/{s_{T+1}})$ and $({\Gamma_T}/{s_{T+1}})$ have the same finite subsequential limits and the
same divergence to $+\infty$. This proves both equalities in
\eqref{eq:key-liminf-limsup}.

Finally,
\[
    \frac{\ln X^{(\pi)}_{T+1}}{s_{T+1}}
    =
    \frac{\ln X_{T+1}^*}{s_{T+1}}
    +
    \frac{
        \ln\pi_{T+1}
        -
        \ln q_{T+1}^{(\delta^*)}
    }{s_{T+1}},
\]
and the last term converges to $\delta^*$ by
\eqref{eq:PQ-first-order}. Combining this with
\eqref{eq:key-liminf-limsup} gives
\eqref{eq:H-liminf-limsup}.
\end{proof}

The following Theorem~\ref{thm:subexp} separates the two sides of the first-order scale. Its lower/upper conditions are equivalent to 
\[ 
\liminf_T\frac{\ln X^{(\pi)}_T}{s_T}\geq\delta^*\quad\text{a.s.}\qquad \text{and}\qquad\limsup_T\frac{\ln X^{(\pi)}_T}{s_T}\leq\delta^*\quad\text{a.s.},
\] 
respectively. Both hold precisely when 
$
\frac{\ln X^{(\pi)}_T}{s_T}\stackrel{a.s.}\longrightarrow\delta^*. 
$ 
This is a first-order statement. Under the stronger second-order product estimate \eqref{eq:PQ-second-order}, exact two-sided control of $X_T^*$ requires an additive comparison:
\[ 
\liminf_T X_T^*=0\;\text{a.s.} \quad\Longleftrightarrow\quad \liminf_T\bigl(\ln X^{(\pi)}_T-\delta^*s_T\bigr)=-\infty\;\text{a.s.},
\] and 
\[ 
\limsup_T X_T^*=+\infty\;\text{a.s.} \quad\Longleftrightarrow\quad \limsup_T\bigl(\ln X^{(\pi)}_T-\delta^*s_T\bigr)=+\infty\;\text{a.s.}.
\] 
The additive comparison needed to obtain these conclusions is what we used in
Section~\ref{sec:main-last}.

\begin{theorem}[First-order lower and upper criteria]\label{thm:subexp} 
Consider the basic recursion~\eqref{eq:recursion}. 
Under the same assumptions as in Lemma~\ref{lem:logZ-identity} (i.e. Assumptions~\ref{ass:reinforcement},  \ref{ass:step-size-1a}, \ref{ass:step-size-1b} and \ref{ass:moment}, 
together with \eqref{eq:log_inc_mean_o}), the
following statements are equivalent on the event $\{X^{(\pi)}_t\to+\infty\}$.

\medskip

\begin{xltabular}{\textwidth}{
    @{}
    >{\raggedright\arraybackslash}X
    >{\raggedright\arraybackslash}X
    @{}
}
\textup{(a) Lower criteria.}
&
\textup{(b) Upper criteria.}
\\[0.75em]

\refstepcounter{lowercriterion}\label{sub-i}%
\textup{\thelowercriterion}
For every $\eps>0$,
\cellmath{
    e^{\eps s_T}X_T^*
    \stackrel{\mathrm{a.s.}}{\longrightarrow}+\infty.
}
&
\refstepcounter{uppercriterion}\label{sub-i-upper}%
\textup{\theuppercriterion}
For every $\eps>0$,
\cellmath{
    e^{-\eps s_T}X_T^*
    \stackrel{\mathrm{a.s.}}{\longrightarrow}0.
}
\\

\refstepcounter{lowercriterion}\label{sub-ii}%
\textup{\thelowercriterion}
For every $\eps>0$,
\cellmath{
    e^{-(\delta^*-\eps)s_T}X^{(\pi)}_T
    \stackrel{\mathrm{a.s.}}{\longrightarrow}+\infty.
}
&
\refstepcounter{uppercriterion}\label{sub-ii-upper}%
\textup{\theuppercriterion}
For every $\eps>0$,
\cellmath{
    e^{-(\delta^*+\eps)s_T}X^{(\pi)}_T
    \stackrel{\mathrm{a.s.}}{\longrightarrow}0.
}
\\

\refstepcounter{lowercriterion}\label{sub-iii}%
\textup{\thelowercriterion}
\cellmath{
    \liminf_T\frac{\Gamma_T}{s_{T+1}}
    \geq 0
    \quad\mathrm{a.s.}
}
&
\refstepcounter{uppercriterion}\label{sub-iii-upper}%
\textup{\theuppercriterion}
\cellmath{
    \limsup_T\frac{\Gamma_T}{s_{T+1}}
    \leq 0
    \quad\mathrm{a.s.}
}
\\

\refstepcounter{lowercriterion}\label{sub-iv}%
\textup{\thelowercriterion}
\cellmath{
    \liminf_T
    \frac{\ln X_{T+1}^*}{s_{T+1}}
    \geq 0
    \quad\mathrm{a.s.}
}
&
\refstepcounter{uppercriterion}\label{sub-iv-upper}%
\textup{\theuppercriterion}
\cellmath{
    \limsup_T
    \frac{\ln X_{T+1}^*}{s_{T+1}}
    \leq 0
    \quad\mathrm{a.s.}
}
\\

\refstepcounter{lowercriterion}\label{sub-v}%
\textup{\thelowercriterion}
\cellmath{
    \liminf_{T\to\infty}
    \frac{\displaystyle\sum_{t\leq T}V_{t+1}/X^{(\pi)}_t}
         {\delta^*s_{T+1}}
    \geq 1     \quad\mathrm{a.s.}
}
&
\refstepcounter{uppercriterion}\label{sub-v-upper}%
\textup{\theuppercriterion}
\cellmath{
    \limsup_{T\to\infty}
    \frac{\displaystyle\sum_{t\leq T}V_{t+1}/X^{(\pi)}_t}
         {\delta^*s_{T+1}}
    \leq 1
        \quad\mathrm{a.s.}
}
\\

\refstepcounter{lowercriterion}\label{sub-vi}%
\textup{\thelowercriterion}
For every $\alpha_1,\alpha_2>0$ with $\alpha_1<\delta^*$,
\cellmath{
    \bigl(X^{(\pi)}_T\bigr)^{-\alpha_2}
    \exp\{(\delta^*-\alpha_1)\alpha_2s_T\}
    \stackrel{a.s.}{\longrightarrow}  0,
}
&
\refstepcounter{uppercriterion}\label{sub-vi-upper}%
\textup{\theuppercriterion}
For every $\alpha_1,\alpha_2>0$,
\cellmath{
    \bigl(X^{(\pi)}_T\bigr)^{\alpha_2}
    \exp\{-(\delta^*+\alpha_1)\alpha_2s_T\}
    \stackrel{a.s.}{\longrightarrow} 0,
}
\end{xltabular}

\medskip

Consequently, on the event $\{X^{(\pi)}_t\to+\infty\}$, the lower and upper conditions hold simultaneously if and
only if
\begin{equation}\label{eq:first-order-exact}
    \frac{\Gamma_T}{s_{T+1}}
    \stackrel{a.s.}\longrightarrow0, \qquad \text{or} \qquad
    \frac{\ln X_{T+1}^*}{s_{T+1}}
    \stackrel{a.s.}\longrightarrow0,
\qquad \text{or} \qquad    
        \frac{
        \sum_{t\leq T}V_{t+1}/X^{(\pi)}_t
    }{
        \delta^*s_{T+1}
    }
    \stackrel{a.s.}\longrightarrow1.
\end{equation}
or, equivalently,
\begin{equation}\label{eq:first-order-exact-H}
    \frac{\ln X^{(\pi)}_T}{s_T}
    \stackrel{a.s.}\longrightarrow\delta^*.
\end{equation}
\end{theorem}

\begin{proof}[Proof of Theorem~\ref{thm:subexp}]
Work on $\{X^{(\pi)}_t\to+\infty\}$. By
Lemma~\ref{lem:logZ-identity}, conditions \ref{sub-iii} and
\ref{sub-iv} are equivalent, as are \ref{sub-iii-upper} and
\ref{sub-iv-upper}. Moreover, \eqref{eq:PQ-first-order} gives
\(
\ln X^{(\pi)}_T=\ln X_T^*+\delta^*s_T+o(s_T). 
\)
Taking logarithms in the displayed exponential conditions (and
absorbing each $o(s_T)$ term by halving the corresponding positive
constant) now gives
\[
\ref{sub-i}\Longleftrightarrow\ref{sub-ii}
\Longleftrightarrow\ref{sub-iv}
\Longleftrightarrow\ref{sub-vi}
\]
and, with the signs reversed,
\[
\ref{sub-i-upper}\Longleftrightarrow\ref{sub-ii-upper}
\Longleftrightarrow\ref{sub-iv-upper}
\Longleftrightarrow\ref{sub-vi-upper}.
\]

It remains to treat the random sums. Set
\[
A_T=\sum_{t\leq T}\frac{\E(V_{t+1}\mid\F_t)}{X^{(\pi)}_t},
\qquad
S_T=\sum_{t\leq T}\frac{V_{t+1}}{X^{(\pi)}_t}.
\]
Theorem~\ref{thm:log-compensator_comp-relative} gives $A_T\stackrel{a.s.}\to+\infty$, and
Proposition~\ref{thm:log-compensator}\ref{comp-replacement}, applied with
$P_t=X^{(\pi)}_t$  gives $S_T/A_T\stackrel{a.s.}\to 1$ a.s. Since
\[
\frac{A_T}{\delta^*s_{T+1}}
=1-\frac{s_\tau}{s_{T+1}}+\frac{\Gamma_T}{\delta^*s_{T+1}}\geq0,
\]
we obtain simultaneously, as \(s_\tau/s_{T+1}\to0\),
\[
\begin{split}
\liminf_T\frac{S_T}{\delta^*s_{T+1}}
&=1+\frac1{\delta^*}
\liminf_T\frac{\Gamma_T}{s_{T+1}}\qquad\text{a.s.},\\
\limsup_T\frac{S_T}{\delta^*s_{T+1}}
&=1+\frac1{\delta^*}
\limsup_T\frac{\Gamma_T}{s_{T+1}}\qquad\text{a.s.}.
\end{split}
\]
Thus \ref{sub-v} is equivalent to \ref{sub-iii}, and
\ref{sub-v-upper} is equivalent to \ref{sub-iii-upper}.

Finally, the lower and upper conditions hold simultaneously precisely
when $\Gamma_T/s_{T+1}\stackrel{a.s.}\to 0$. The identities above and
Lemma~\ref{lem:logZ-identity} show that this is equivalent to each
condition in \eqref{eq:first-order-exact} and to
\eqref{eq:first-order-exact-H}.
\end{proof}
The next proposition gives symmetric sufficient conditions for the lower and upper parts of Theorem~\ref{thm:subexp}. An asymptotically negligible negative relative perturbation yields the lower first-order criteria, while an asymptotically negligible positive relative perturbation yields the corresponding upper criteria.

\begin{proposition}[First-order criteria from negligible relative perturbations]\label{prop:relative-negative-subexp}
Consider the basic recursion~\eqref{eq:recursion}. 
Under the same assumptions as in Theorem~\ref{thm:subexp}, on the event
$\{X^{(\pi)}_t\to +\infty\}$, if
\(
    {\rho_t^-}/{X_t}\stackrel{a.s.}\longrightarrow 0,
\)
(resp.,
\(
    {\rho_t^+}/{X_t}\stackrel{a.s.}\longrightarrow 0,
\))
then the equivalent lower conditions
\ref{sub-i}--\ref{sub-vi}
(resp.\ the equivalent upper conditions \ref{sub-i-upper}--\ref{sub-vi-upper})
of Theorem~\ref{thm:subexp} hold true.
\end{proposition}
\begin{proof}[Proof of Proposition~\ref{prop:relative-negative-subexp}]
Work on $\{X^{(\pi)}_t\to+\infty\}$.
Set
\(
\eps_t^-:={\rho_t^-}/{X_t}
\) (resp. $\eps_t^+:={\rho_t^+}/{X_t}$). 
By hypothesis, $\eps_t^-\stackrel{a.s.}\to 0$ (resp.\ $\eps_t^+\stackrel{a.s.}\to 0$), and the definition of $\Gamma_T$
gives
\[
\Gamma_T
\geq-\sum_{t\leq T}\htr_{t+1}\eps_t^-
\qquad
\Bigl(\text{resp. }\Gamma_T
\leq\sum_{t\leq T}\htr_{t+1}\eps_t^+\Bigr).
\]
The weighted Cesaro theorem and $s_{T+1}\to+\infty$ yield
\[
\frac1{s_{T+1}}
\sum_{t\leq T}\htr_{t+1}\eps_t^-
\longrightarrow0
\qquad
\Bigl(\text{resp. }\frac1{s_{T+1}}
\sum_{t\leq T}\htr_{t+1}\eps_t^+
\longrightarrow 0\Bigr)
\]
Consequently,
\[
\liminf_T\frac{\Gamma_T}{s_{T+1}}\geq 0\quad\text{a.s.}
\qquad
\Bigl(\text{resp. }\limsup_T\frac{\Gamma_T}{s_{T+1}}\leq 0\quad\text{a.s.}\Bigr).
\]
This is condition~\ref{sub-iii} (resp.\
condition~\ref{sub-iii-upper}) of Theorem~\ref{thm:subexp}, which
proves all the equivalent lower conditions
\ref{sub-i}--\ref{sub-vi} (resp. upper conditions
\ref{sub-i-upper}--\ref{sub-vi-upper}).
\end{proof}

For the convergent exact-normalization theorem, a convenient condition controlling the negative part of the perturbation is 
\begin{equation}\label{eq:rho-negative} 
\sum_t \htr_{t+1}\frac{\rho_t^-}{X_t}<+\infty \qquad\text{a.s. on }\{X^{(\pi)}_t\to+\infty\}. 
\end{equation} 
Together with the corresponding positive control, this leads to the absolute relative summability appearing in Theorem~\ref{thm:positive-all_second}.
However, it remains useful to have a condition on the negative perturbation that does not involve the random denominator $X_t$. The next proposition provides such a criterion. It is independent of Assumption~\ref{ass:quadratic-relative-perturbation}, which controls the positive predictable second-order contribution in Theorem~\ref{thm:main-positive-all}.

\begin{proposition}[Criterion for controlling the negative perturbation]\label{prop:rho-negative-sufficient} 
Consider the basic recursion~\eqref{eq:recursion}. 
Under the same assumptions as in Lemma~\ref{lem:logZ-identity}, assume also
\begin{equation}\label{eq:rho-negative-sufficient}
        \rho_t^-=O(\pi_t^{-1})
        \qquad\mbox{on }\{X^{(\pi)}_t\to+\infty\}.
\end{equation}
Then condition~\eqref{eq:rho-negative} holds true.
\end{proposition}

\begin{proof}[Proof of Proposition~\ref{prop:rho-negative-sufficient}]
Work on $\{X^{(\pi)}_t\to+\infty\}$.  Since $X^{(\pi)}_t=\pi_tX_t$, condition
\eqref{eq:rho-negative-sufficient} gives
\[
        \frac{\rho_t^-}{X_t}=O(\bigl(X^{(\pi)}_t\bigr)^{-1})=o(1).
\]
Therefore, by Proposition~\ref{prop:relative-negative-subexp}, the equivalent lower conditions of Theorem~\ref{thm:subexp} hold. Taking, for instance, $\alpha_1=\delta^*/2$ and $\alpha_2=1$ in condition~\ref{sub-vi} gives 
\[ 
\bigl(X^{(\pi)}_t\bigr)^{-1} =o\bigl(e^{-\delta^*s_t/2}\bigr) =o\bigl(e^{-\delta^*s_{t+1}/2}\bigr). 
\]

Then, using again \eqref{eq:rho-negative-sufficient},
and \eqref{eq:exp-sum} with
$\eps=\delta^*/2$, we obtain \eqref{eq:rho-negative}, since
\[
        \sum_{t\ge0}\htr_{t+1}\frac{\rho_t^-}{X_t} 
        =O\Bigl(\sum_{t\ge0}\frac{\htr_{t+1}}{X^{(\pi)}_t}\Bigr) 
        =O\Bigl(\sum_{t\ge0}\htr_{t+1} e^{-\delta^*s_{t+1}/2} \Bigr) 
        <+\infty\quad\text{a.s.}\,. \qedhere
\]
\end{proof}

\section{Examples and Counterexamples}\label{sec:rates} 

\subsection{The behavior of the normalization $q_t^{(\delta^*)}$}\label{subsec:Q-regimes}
Theorem~\ref{thm:positive-all_second-maringale-noise} affirms, under suitable assumptions, that 
 $X_t\stackrel{a.s.}\sim X_\infty^*/q_t^{(\delta^*)}$.  Hence, for the sequence $(X_t)$ we can have 
 decay to zero, convergence to a positive finite limit, growth to infinity, or even non-convergence. 
 Remark~\ref{cor:polynomial} provides examples for the first three regimes. The following example shows the non-convergence regime. 

\begin{example}[A non-convergent normalization]\label{ex:Q-oscillating}
For the basic recursion~\eqref{eq:recursion}, set $Y_{t+1}:=\delta^*X_t$ (which trivially implies $\rho_t\equiv 0$),
\[
        r_t:=\frac{b}{t+\theta},
        \qquad
        p_t:=\frac{b\delta^*+\kappa\sigma_t}{t+\theta},
        \qquad 0<\kappa<b\delta^*, \qquad \theta>b\delta^*+\kappa,
\]
where $\sigma_t=(-1)^m$ for $2^m\le t<2^{m+1}$. Then $p_t-\delta^*r_t={\kappa\sigma_t}/{(t+\theta)}$, 
whose partial sums have two distinct sub-sequential limits and so   
$q_t^{(\delta^*)}$ does not converge. Moreover, we have $
        X^*_t=q_{t}^{(\delta^*)}X_{t}
        =q_0^{(\delta^*)}X_0$ 
and, by \eqref{eq:PQ-first-order}, $\ln(\pi_t/q_t^{(\delta^*)})\sim\delta^*s_t\to+\infty$, so
$X^{(\pi)}_t=\pi_t X_t=(\pi_t/q_t^{(\delta^*)})X^*_t\to +\infty$ a.s.\,.  Thus, when $X_0>0$,  the normalized process $X^*_t$ 
has a strictly positive limit although neither $q_t^{(\delta^*)}$ nor $X_t$ converges.
\end{example}

\subsection{Counterexamples for Theorem~\ref{thm:main-positive-all-maringale-noise} when canceling an assumption}
The first deterministic example shows that 
one may have the correct first-order behavior and even $\Gamma_T\to 0$,
while $X_T^*\to 0$. More precisely, in the following example 
all the assumptions of Theorem~\ref{thm:main-positive-all-maringale-noise} except
Assumption~\ref{ass:quadratic-relative-perturbation} hold, and we get 
\[
    \Gamma_T\to 0,
    \qquad
    X^{(\pi)}_t\to+\infty,
    \qquad
    X_t^*\to 0.
\]
The obstruction is the accumulation of second-order
logarithmic losses. 

\begin{example}[A counterexample for Theorem~\ref{thm:main-positive-all-maringale-noise} if Assumption~\ref{ass:quadratic-relative-perturbation} does not hold]
\label{ex:deterministic-counterexample}
Define
\[
    z_0>0,\qquad z_{t+1}=z_t(1+u_t)\,,
\]
where we choose $(u_t)$ in disjoint blocks. Specifically, fix $c>1$ and
$\beta\in (0,1)$, and put $a_k=c/\sqrt{k}$. Choose
$n_k\geq e^k$ with disjoint successive blocks
$I_k=\{n_k,\ldots,m_k\}$. Set
\[
    u_{n_k}:=a_k,
    \qquad
    u_t:=-\frac{\beta}{t+2},
    \quad n_k<t<m_k,
\]
where $m_k$ is the first index for which the negative weighted
increments can compensate the positive one, and choose
$u_{m_k}\in[-\beta/(m_k+2),0]$ so that
\begin{equation}\label{eq:counter-cancellation}
    \sum_{t=n_k}^{m_k}
    \frac{t+2}{t+1}u_t=0.
\end{equation}
Set $u_t:=0$ outside the blocks. Note that $u_t\to 0$. 
Now, let
\[
    p_{t+1}:=\frac{2}{t+3},
    \qquad
    r_{t+1}:=\frac{1}{t+3}.
\]
Moreover, define 
$$
X_t:=x_t:=\frac{2z_t}{t+2}
$$ 
and, taking $\delta^*:=1$, $\rho_t:=2z_t u_t$ and so $Y_{t+1}:=y_{t+1}:=x_t+\rho_t=2z_t(1/{(t+2)}+u_t)$, we obtain that 
$(x_t)_t$ satisfies the basic recursion~\eqref{eq:recursion} with $y_t \geq 0$ (because $u_t\geq-\beta/(t+2)$) 
and also Assumption~\ref{ass:reinforcement} is satisfied.  We have 
\[
    \pi_t=\frac{(t+1)(t+2)}2,
    \qquad
    q_t^{(1)}=\frac{t+2}{2},
    \qquad
    \htr_{t+1}=\frac1{t+1},
\]
and so Assumption~\ref{ass:step-size} holds. Since the model is deterministic,
Assumption~\ref{ass:moment} is automatic. Positive perturbations occur only at $n_k$,
and
\[
    q_{n_k+1}^{(1)}r_{n_k+1}\rho_{n_k}^+
    =
    z_{n_k}a_k
    \asymp
    k^{-(c^2+1)/2}.
\]
Since $c>1$, Assumption~\ref{ass:lim_sup2} holds. Moreover, as $X^{(\pi)}_t=h_t=(t+1)z_t$, we have  
\[
    h_{n_k}
    =(n_k+1)z_{n_k}\longrightarrow+\infty,
\]
and hence $X^{(\pi)}_t=h_t\to+\infty$. We also have 
\[    
\frac{\pi_{t+1}r_{t+1}\rho_t}{h_t}
    =
    \frac{t+2}{t+1}u_t.
\]
Hence \eqref{eq:counter-cancellation} makes the total contribution of
each block to $\Gamma_T$ equal to zero; therefore
\[
    \Gamma_T\longrightarrow 0.
\]
On the other hand, \eqref{eq:counter-cancellation} gives
\[
    \sum_{t\in I_k}u_t
    =
    O\left(\frac{a_k}{n_k}\right),
    \qquad
    \sum_{\substack{t\in I_k\\t\neq n_k}}u_t^2
    =
    O\left(\frac{a_k}{n_k}\right),
\]
while $\sum_{t\in I_k}|u_t|=O(a_k)$. Therefore,
\[
\begin{aligned}
    \sum_{t\in I_k}\ln(1+u_t)
    =
    -\frac12a_k^2
    +O\left(a_k^3+\frac{a_k}{n_k}\right)  =
    -\frac{c^2}{2k}+O(k^{-3/2}).
\end{aligned}
\]
Summing over the blocks yields
\[
    \ln z_{n_k}
    =
    -\frac{c^2}{2}\ln k+O(1),
    \qquad
    z_{n_k}\asymp k^{-c^2/2}.
\]
Since the maximal increase inside the $k$-th block is only the factor
$1+a_k\to1$, it follows that
\[
    X_t^*=x_t^*=q_t^{(1)}x_t=z_t\longrightarrow 0.
\]
Note also that 
\[
d_t=\frac{\E(V_{t+1}\mid\F_t)}{h_t}=\frac1{t+1}+\frac{t+2}{t+1}u_t\to0,
\]
so \eqref{eq:log_inc_mean_o} holds.
\\
\indent It remains to check that the failure is precisely that of
Assumption~\ref{ass:quadratic-relative-perturbation}. At $t=n_k$,
\[
    \htr_{n_k+1}
    \frac{\rho_{n_k}^+}{X_{n_k}}
    =
    \frac{n_k+2}{n_k+1}a_k
    \sim\frac{c}{\sqrt{k}},
\]
and therefore Assumption~\ref{ass:quadratic-relative-perturbation} fails. The issue is 
the accumulated quadratic logarithmic remainder.
Indeed, we have 
\[
    \sum_{t\leq T}d_t
    -\bigl(\ln h_{T+1}-\ln h_0\bigr)
    =
    \sum_{t\leq T}\bigl[d_t-\ln(1+d_t)\bigr].
\]
At $t=n_k$, we have $d_{n_k}\sim a_k$, and hence
\[
    d_{n_k}-\ln(1+d_{n_k})
    \sim\frac{a_k^2}{2}
    =
    \frac{c^2}{2k}
\]
and so the corresponding series diverges. 
\end{example}

\begin{remark}[Different roles for the second-order assumptions]
The two second-order assumptions, i.e.~Assumption~\ref{ass:moment} and Assumption~\ref{ass:quadratic-relative-perturbation}, 
control different parts of the logarithmic argument.  
Assumption~\ref{ass:moment} controls the martingale fluctuations of $\ln X^{(\pi)}_t$.
Assumption~\ref{ass:quadratic-relative-perturbation}, together with
Assumption~\ref{ass:step-size-2},
controls the predictable quadratic Taylor term because
\[
        0\le
        \frac{\E(V_{t+1}\mid\F_t)}{X^{(\pi)}_t}
        \le
        \delta^*\htr_{t+1}
        +\htr_{t+1}\frac{\rho_t^+}{X_t}.
\]
 The counterexample in
Example~\ref{ex:deterministic-counterexample} has zero conditional variance (and so Assumption~\ref{ass:moment} trivially holds true), but
violates Assumption~\ref{ass:quadratic-relative-perturbation}. Hence, the accumulated predictable quadratic loss forces $X_t^*\to 0$.  
Thus Assumption~\ref{ass:moment} and
Assumption~\ref{ass:quadratic-relative-perturbation} address genuinely distinct mechanisms.
\end{remark}

\begin{remark}
Example~\ref{ex:deterministic-counterexample} can be slightly modified so that all the hypotheses of Theorem~\ref{thm:main-positive-all}, except Assumption~\ref{ass:quadratic-relative-perturbation}, still hold, while
\[
    \limsup_{T\to\infty}\Gamma_T=+\infty,
    \qquad
    \liminf_{T\to\infty}\Gamma_T=-\infty,
    \qquad
    X_T^*\longrightarrow0.
\]
The idea is to replace each compensating block by alternating positive and negative blocks whose total contributions to $\Gamma_T$ have increasing magnitude, while keeping the individual increments small and accumulating enough quadratic logarithmic loss to force $X_T^*\to0$. Thus, without Assumption~\ref{ass:quadratic-relative-perturbation}, even two-sided unbounded oscillations of the compensator $\Gamma_T$ are compatible with the convergence to zero of the normalized process $X_T^*$.
\end{remark}

The second deterministic example shows that  $X^*_\infty>0$ does not hold if condition~\ref{pos-VI2} of Theorem~\ref{thm:positive-all_second-maringale-noise} is not satisfied, 
even if 
$\rho_t^-=O(t^\eps/\pi_t)$ for a small $\eps>0$ (hence the bound in Proposition~\ref{prop:rho-negative-sufficient} is sharp). 

\begin{example}[A counterexample for Theorem~\ref{thm:positive-all_second-maringale-noise} if condition~\ref{pos-VI2} 
does not hold]\label{ex:counterexample_limit_0}
 Fix
$\eps\in(0,1/2)$ and set
\[
        \delta^*:=1,
        \qquad b:=2\eps,
        \qquad \lambda:=1-2\eps,
        \qquad \eta:=1-\eps.
\]
Let
\begin{equation}\label{eq:sharp-rates}
        p_{t+1}:=\frac1{t+2},
        \qquad
        r_{t+1}:=\frac{b}{t+2}.
\end{equation}
Then
        $\pi_t\asymp t$ and 
        $q_t^{(\delta^*)}\asymp t^\lambda$.
Choose $c>0$ and, for all large $t$, set
$\rho_t:=-c(t+2)^{-\eta}$.  Define
\begin{equation}\label{eq:sharp-Z-tail}
        X_t^*:=x_t^*:= q_t^{ ( \delta^*) } x_t
        =c\sum_{k=t}^{\infty}
        q_{k+1}^{(\delta^*)}r_{k+1}(k+2)^{-\eta}.
\end{equation}
Since the summand is of order $k^{-1-\eps}$, the series is finite and
$X_t^*=x^*_t\to 0$ (hence condition~\ref{pos-I2} of Theorem~\ref{thm:positive-all_second-maringale-noise} fails).  The identity
$q_t^{(\delta^*)}=(1-p_{t+1}+r_{t+1})q_{t+1}^{(\delta^*)}$ shows that the
corresponding recursion has
$Y_{t+1}:=y_{t+1}:=x_t-c(t+2)^{-\eta}$.
Note that 
\[
        x_t\sim\frac{bc}{\eta-\lambda}t^{-\eta}=2ct^{-\eta},
\]
so $Y_{t+1}\ge0$ eventually; finitely many initial terms can be adjusted.
The example is deterministic, hence Assumption~\ref{ass:moment} holds, while
$\rho_t^+=0$ makes Assumptions~\ref{ass:quadratic-relative-perturbation} and \ref{ass:lim_sup2} automatic.  Moreover,
$ X^{(\pi)}_t=h_t\asymp t^\eps\to +\infty$. 
Finally, $\pi_t^{-1}\asymp t^{-1}$ and $\rho_t^-\asymp t^{-1+\eps}$,  
so \eqref{eq:rho-negative-sufficient} fails by the factor $t^\eps$.
More importantly,
\[
        \htr_{t+1}\frac{\rho_t^-}{X_t}\asymp t^{-1},
\]
and therefore \eqref{eq:rho-negative} fails.  Equivalently, condition~\ref{pos-VI2}
of Theorem~\ref{thm:positive-all_second-maringale-noise} fails. As a consequence, in this case, Theorem~\ref{thm:positive-all_second-maringale-noise}\ref{pos-V2} and Theorem~\ref{thm:main-positive-all-maringale-noise} imply $\lim_T \Gamma_T=-\infty$.
\end{example}

\subsection{Example for $(X^{(\pi)}_t)$ converging to zero}\label{sec:H_t_converge_zero}
According to Theorem~\ref{cor-X0}(b), for the modified recursion~\eqref{eq:recursion:mart} we have a dichotomy: 
$X^{(\pi)}_t\to 0$ or $X^{(\pi)}_t\to +\infty$. The next example we shows a case in which the a.s. convergence of $X^{(\pi)}_t$ toward zero. 
\begin{example}
Set
\(
p_{t+1}:=0,\  r_{t+1}:=\frac1{t+2},\
X_0:=1,\  Y_{t+1}:=X_t,\ 
\delta^*:=1,\ \rho_t:=0.
\)
Independently at each time, choose the multiplier
\[
R_{t+1}:=
\begin{cases}
\dfrac{(t+3)^2}{2},
   &\text{with probability } r_{t+1}^{2},\\[2mm]
\dfrac12,
   &\text{with probability }1-r_{t+1}^{2},
\end{cases}
\qquad X_{t+1}:=R_{t+1}X_t.
\]
Since \(\E(R_{t+1}\mid\F_t)=1+r_{t+1}\), setting
\[
\Delta N_{t+1}
 :=\frac{R_{t+1}-1-r_{t+1}}{r_{t+1}}
\]
gives a martingale difference and puts the update in the form
\eqref{eq:recursion:mart}. Moreover, \(X^{(\pi)}_t=X_t>0\),
\(\htr_{t+1}=r_{t+1}\), and
\(
\E(\widetilde V_{t+1}\mid\F_t)=r_{t+1}X^{(\pi)}_t.
\)
Thus Assumptions~\ref{ass:reinforcement} and
\ref{ass:step-size-1b} hold, and the event \(E\) in
\eqref{eq-dicotomy-H} is almost sure. The probabilities of the large multipliers are summable: indeed, 
\(\sum_t r_{t+1}^2<+\infty\). By Borel--Cantelli, only
finitely many occur almost surely. Thereafter \(X^{(\pi)}_t\) is halved
at every step, so \(X^{(\pi)}_t\stackrel{a.s.}\to 0\). Finally,
\[
\frac{\Var(\widetilde V_{t+1}\mid\F_t)}
     {\E(\widetilde V_{t+1}\mid\F_t)}
=
\frac{(1-r_{t+1}^{\,2})(1+2r_{t+1})^2}
     {4r_{t+1}^{\,3}}\,X^{(\pi)}_t
\stackrel{a.s.}\longrightarrow 0.
\]
Indeed, after the last large multiplier, \(X^{(\pi)}_t\) decays
geometrically, whereas \(r_{t+1}^{-3}\) grows only
polynomially. Thus, \eqref{eq:signed-variance-drift}
holds.
\end{example}

\section{Applications}\label{sec:applications}
We are going to show some applications of the previous theoretical results. It is worthwhile to note that some of the 
considered models are new extensions of known models (e.g. triangular urns with sampling bias and generalized Simon innovation model). 
The applications regard both  
basic recursions~\eqref{eq:recursion} and modified recursions with a multiplicative 
martingale noise~\eqref{eq:recursion:mart}.
In each case the proof reduces to identifying the scalar rate, 
establishing $X^{(\pi)}_t\stackrel{a.s.}\to+\infty$, and checking the relevant noise and
perturbation series.  When $\rho_t=0$, Assumption~\ref{ass:lim_sup2} and
condition~Theorem~\ref{thm:positive-all_second-maringale-noise}\ref{pos-VI2} are automatic. 
When $\rho_t\geq 0$, the lower bound on $\Gamma_T$ is automatic, and condition~Theorem~\ref{thm:positive-all_second-maringale-noise}\ref{pos-V2} reduces the problem to verify 
Assumptions~\ref{ass:quadratic-relative-perturbation} and~\ref{ass:lim_sup2}. 
For the basic recursions, we verify Assumption~\ref{ass:moment} before applying Theorem~\ref{cor-X0}\ref{cor-X0:a}; for modified recursions with martingale noise, we use part~\ref{cor-X0:b} of that theorem after checking~\eqref{eq:signed-variance-drift}.

\subsection{A linear reward-inaction learning automaton}\label{subsec:lri}
The linear reward-inaction rule is a classical learning-automaton
scheme \cite{NarendraThathachar1989}. In particular, the two-armed bandit model
and its convergence and rate theory has been widely studied (e.g. see~\cite{LambertonPages2008,LambertonPagesTarres2004}).
Consider a two-action stochastic environment.  Action $1$ is inferior and
Action $2$ is optimal.  If Action $i$ is selected, it is rewarded with
probability $c_i$, with $0<c_1<c_2<1$. Set $\Deltac:=c_2-c_1$.  
Let $X_t$ be the probability assigned to Action 1.  Under the linear
reward-inaction rule with gain $\gamma_{t+1}$,
\begin{equation}\label{eq:lri-update}
        X_{t+1}=(1-\gamma_{t+1})X_t+\gamma_{t+1}Y_{t+1},
\end{equation}
where
\[
Y_{t+1}:=
\begin{cases}
1, & \text{if action $1$ is selected and rewarded},\\
0, & \text{if action $2$ is selected and rewarded},\\
X_t, & \text{if no reward occurs}.
\end{cases}
\]
Thus $(X_t)$ satisfies dynamics~\eqref{eq:recursion} with $p_{t+1}=r_{t+1}=\gamma_{t+1}$ and
\begin{align}\label{eq:lri-conditional-mean}
        \E(Y_{t+1}\mid\F_t)
        &=(1-\Deltac)X_t+\Deltac X_t^2
        =\delta^*X_t+\rho_t,
\end{align}
where $\delta^*:=1-\Deltac$ and  $\rho_t:=\Deltac X_t^2\ge0$. Hence, Assumption~\ref{ass:reinforcement} is satisfied. 

This harmonic-gain case illustrates our criteria while recovering
the positive slow-rate limit studied in~\cite{LambertonPages2008}.
\begin{proposition}\label{prop:lri}
Suppose $X_0\in(0,1)$ and
\[
        \gamma_{t+1}=\frac{a}{t+\theta},
        \qquad 0<a\le1,
\]
where $\theta>0$ is such that $\gamma_{t+1}\in(0,1)$ for each $t$.  Then, there exists a strictly positive real random variable $\widetilde{X}_\infty$ such that 
\begin{equation}\label{eq:lri-rate}
        t^{a\Deltac}X_t
        \stackrel{a.s.}\longrightarrow
        \widetilde{X}_\infty.
\end{equation}
\end{proposition}

\begin{proof}
We have
\[
\htr_{t+1}
=
\frac{\gamma_{t+1}}{1-\gamma_{t+1}}
=
\frac{a}{t+\theta-a}
\]
so Assumptions~\ref{ass:step-size} hold. Since
\(0\leq Y_{t+1}\leq1\) and $
\pi_t=O(t^a)$, Remark~\ref{ex:bounded} gives Assumption~\ref{ass:moment}.
Moreover, $X_0>0$ and \(\rho_t=\Deltac X_t^2\geq0\), so Theorem~\ref{cor-X0}\ref{cor-X0:a} gives \(X^{(\pi)}_t\to+\infty\) a.s.\,. 
Now, the first step is to prove that
\begin{equation}\label{eq:lri-X-zero-direct}
X_t\stackrel{a.s.}\longrightarrow0.
\end{equation}
To this purpose, note that 
\[
\E(X_{t+1}-X_t\mid\F_t)
=
-\Deltac\gamma_{t+1}X_t(1-X_t)\leq0.
\]
Thus \((X_t)\) is a bounded supermartingale and so it converges almost surely
to a random variable \(X_\infty\in[0,1]\). Moreover,
\[
\sum_t\gamma_{t+1}X_t(1-X_t)<+\infty
\qquad\text{a.s.}
\]
Since \(\sum_t\gamma_{t+1}=+\infty\), this implies 
$X_\infty\in\{0,1\}$ a.s.\,.
It remains to exclude \(X_\infty=1\). Set
$Z_t:=1-X_t$
 and 
$Y^Z_{t+1}:=1-Y_{t+1}$. 
Then \(Z_0>0\), \(0\leq Y^Z_{t+1}\leq1\), and
\[
Z_{t+1}
=
(1-\gamma_{t+1})Z_t
+\gamma_{t+1}Y^Z_{t+1}.
\]
Furthermore,
\[
\begin{aligned}
\E(Y^Z_{t+1}\mid\F_t)
=
1-\E(Y_{t+1}\mid\F_t)
=
(1+\Deltac)Z_t-\Deltac Z_t^2.
\end{aligned}
\]
Thus \((Z_t)\) follows the basic recursion~\eqref{eq:recursion} with 
$\delta_Z^*:=1+\Deltac$ and 
$\rho_t^Z:=-\Deltac Z_t^2$. 
In particular,
\[
\frac{(\rho_t^Z)^-}{Z_t}
=
\Deltac Z_t
\leq\Deltac
<
1+\Deltac
=
\delta_Z^*.
\]
Since \(0\leq Y^Z_{t+1}\leq1\) and \(\pi_{t+1}r_{t+1}=O(1)\), Remark~\ref{ex:bounded} also gives Assumption~\ref{ass:moment} for the recursion of \(Z_t\). Theorem~\ref{cor-X0}\ref{cor-X0:a} therefore gives $Z_t^{(\pi)}:=\pi_tZ_t\stackrel{a.s.}\longrightarrow+\infty$. 
Moreover, since \(0<Z_t\leq1\),
\[
0\leq
\frac{\E(V^Z_{t+1}\mid\F_t)}{Z_t^{(\pi)}}
=\htr_{t+1}(1+\Deltac-\Deltac Z_t)
\leq(1+\Deltac)\htr_{t+1}\longrightarrow0.
\]
Thus condition~\eqref{eq:log_inc_mean_o} holds for \(Z\). 
Consider now the event \(\{X_\infty=1\}\). On this event,
\(Z_t\stackrel{a.s.}\to0\), and hence
\[
\frac{(\rho_t^Z)^-}{Z_t}
=
\Deltac Z_t\longrightarrow0,
\qquad
\frac{(\rho_t^Z)^+}{Z_t}=0.
\]
Proposition~\ref{prop:relative-negative-subexp} therefore gives both
the lower and upper criteria of Theorem~\ref{thm:subexp}. Consequently,
\[
\frac{\ln Z_t^{(\pi)}}{s_t}
\longrightarrow
\delta_Z^*
=
1+\Deltac
\qquad\text{on }\{X_\infty=1\}.
\]
On the other hand,
\[
\ln\pi_t
=
-\sum_{k<t}\ln(1-\gamma_{k+1}),
\qquad
s_t
=
\sum_{k<t}\frac{r_{k+1}}{1-p_{k+1}}=
\sum_{k<t}\frac{\gamma_{k+1}}{1-\gamma_{k+1}}.
\]
Since \(\sum_t\gamma_{t+1}^2<+\infty\), we have 
\(
{\ln\pi_t}/{s_t}\rightarrow 1.
\)
Therefore,
\[
\frac{\ln Z_t}{s_t}
=
\frac{\ln Z_t^{(\pi)}-\ln\pi_t}{s_t}
\longrightarrow
\Deltac>0
\qquad\text{on }\{X_\infty=1\}.
\]
This is impossible because \(0<Z_t\leq1\), and hence
\(\ln Z_t\leq0\). Thus 
$\Pbb(X_\infty=1)=0$, which proves \eqref{eq:lri-X-zero-direct}. 
 We are now ready to verify that Assumption~\ref{ass:lim_sup2} holds. Fix
\(\eps\in(0,\Deltac/2)\), put
\(\delta_\eps:=\delta^*+\eps\), and rewrite
\[
\E(Y_{t+1}\mid\F_t)
=
\delta_\eps X_t+\rho_t^{(\eps)},
\qquad
\rho_t^{(\eps)}
=
X_t(\Deltac X_t-\eps).
\]
Since \(X_t\stackrel{a.s.}\to 0\), one has
\(\rho_t^{(\eps)}\leq0\) eventually almost surely. Thus
Assumption~\ref{ass:lim_sup2} holds for this decomposition, and
Theorem~\ref{thm:finite-limit} gives
$q_t^{(\delta_\eps)}X_t=O(1)$. 
Since, by Remark~\ref{cor:polynomial},  
\(
{q_t^{(\delta_\eps)}}
/ {t^{a(\Deltac-\eps)}}
\longrightarrow \ell_\eps > 0
\),
we obtain 
$X_t=O\bigl(t^{-a(\Deltac-\eps)}\bigr)$. Returning to the original decomposition, we get 
\[
\sum_t
\htr_{t+1}\frac{|\rho_t|}{X_t}
=
\Deltac\sum_t\htr_{t+1}X_t
<+\infty\quad\mbox{a.s.}\,.
\]
Theorem~\ref{prop:absolute-relative-criterion}, together with Remark~\ref{cor:polynomial}, therefore gives the thesis because  
\(1-\delta^*=\Deltac\).
\end{proof}

\subsection{Triangular urn with a signed sampling bias}
\label{subsec:triangular-urn}\label{subsec:perturbed-urn}

Fix \(s>0\) and \(a\in(0,s)\), and define the deterministic total mass $\tau_t:=\tau_0+st$. 
Let $X_t:={W_t}/{\tau_t}$, where \(W_t\) is the white mass and \(0<W_0\leq\tau_0\).
We introduce a signed perturbation $\eta_t$ in the sampling
probability. Let \(I_{t+1}\) be the indicator that a white ball is drawn.
Conditionally on \(\F_t\), assume that \(I_{t+1}\) is Bernoulli with
parameter
\[
\Pbb(I_{t+1}=1\mid\F_t)=X_t+\eta_t,
\]
where \((\eta_t)\) is $\F_t$-measurable, possibly signed, and satisfies
\begin{equation}\label{eq:triangular-bias-range}
-X_t\leq\eta_t\leq1-X_t.
\end{equation}
If a white ball is drawn, \(a\) units of white mass and \(s-a\) units
of black mass are added; otherwise, \(s\) units of black mass are
added. Thus
\[
W_{t+1}=W_t+aI_{t+1},
\]
and
\begin{equation}\label{eq:triangular-recursion}
X_{t+1}
=
\left(1-\frac{s}{\tau_{t+1}}\right)X_t
+\frac{a}{\tau_{t+1}}I_{t+1}.
\end{equation}
The unperturbed triangular urn is recovered by taking
\(\eta_t=0\) \cite{Janson2026}. 
The following proposition extends the strict positive limit for the unperturbed
triangular urn~\cite{Janson2026} to signed predictable sampling biases.
\begin{proposition} 
\label{prop:triangular-signed-bias}
Suppose that
\begin{equation}\label{eq:bias-summability}
\sum_{t}
q_{t+1}^{(1)}
\frac{a}{\tau_{t+1}}\eta_t^+
<+\infty
\qquad\text{a.s.},
\end{equation}
\begin{equation}\label{eq:bias-square-summability}
(a)\quad \sum_{t}(\eta_t^+)^2<+\infty
\qquad\text{a.s.}
\quad 
\text{and}\quad  
(b)\quad\sum_{t}\frac{\eta_t^-}{W_t}<+\infty
\qquad\text{a.s.}
\end{equation}
Then there exist two strictly positive real random variables $
\widetilde{X}_\infty,\, \widetilde{W}_\infty$ such that
\begin{equation}\label{eq:triangular-rate}
t^{(s-a)/s}X_t
\stackrel{a.s.}\longrightarrow
\widetilde{X}_\infty
\end{equation}
and
\begin{equation}\label{eq:triangular-W-rate}
t^{-a/s}W_t
\stackrel{a.s.}\longrightarrow
\widetilde{W}_\infty.
\end{equation}
\end{proposition}

\begin{proof} Consider~\eqref{eq:triangular-recursion} and 
set
\[
p_{t+1}:=\frac{s}{\tau_{t+1}},
\qquad
r_{t+1}:=\frac{a}{\tau_{t+1}},
\qquad
Y_{t+1}:=I_{t+1},
\qquad
\delta^*:=1,
\qquad
\rho_t:=\eta_t.
\]
Then
\[
1-p_{t+1}=\frac{\tau_t}{\tau_{t+1}},
\qquad
\pi_t=\frac{\tau_t}{\tau_0},
\qquad
X^{(\pi)}_t=\frac{W_t}{\tau_0},
\qquad
\htr_{t+1}=\frac{a}{\tau_t}.
\]
Thus Assumptions~\ref{ass:step-size} hold. Since
\(0\leq Y_{t+1}\leq1\) and $\pi_{t+1}r_{t+1}={a}/{\tau_0}$, 
Assumption~\ref{ass:moment} follows from Remark~\ref{ex:bounded}.

We first prove that \(W_t\to+\infty\) using Theorem~\ref{cor-X0}\ref{cor-X0:a}. Indeed, since \(W_0>0\),  
Assumption~\ref{ass:moment} holds, and
\[
\begin{aligned}
\sum_t \htr_{t+1}
 \left(1+\frac{\rho_t}{X_t}\right)
=
a\sum_t\frac1{\tau_t}
+a\sum_t\frac{\eta_t}{W_t}
\geq
a\sum_t\frac1{\tau_t}
-a\sum_t\frac{\eta_t^-}{W_t}
=+\infty
\qquad\text{a.s.}
\end{aligned}
\]
(here the last step uses
\eqref{eq:bias-square-summability}\textup{(b)}), 
Theorem~\ref{cor-X0}\ref{cor-X0:a} gives
\(X^{(\pi)}_t=W_t/\tau_0 \stackrel{a.s.}\to+\infty\).
 Condition~\eqref{eq:bias-summability} is precisely
Assumption~\ref{ass:lim_sup2}. Furthermore,
\[
\htr_{t+1}\frac{\rho_t^+}{X_t}
=
\frac{a\eta_t^+}{W_t}
\leq
\frac{a}{W_0}\eta_t^+,
\]
so \eqref{eq:bias-square-summability}$(a)$ gives
Assumption~\ref{ass:quadratic-relative-perturbation}. Finally, by~\eqref{eq:bias-square-summability}$(b)$, 
\[
\Gamma_T
=
a\sum_{t\leq T}\frac{\eta_t}{W_t},
\qquad
\sup_T \Gamma_T^-
\leq
a\sum_t\frac{\eta_t^-}{W_t}<+\infty.
\]
Thus condition~\ref{pos-V2} of Theorem~\ref{thm:positive-all_second-maringale-noise} holds. Theorem
\ref{thm:positive-all_second} and Remark~\ref{cor:polynomial}, together with~\(\tau_t/t\to s\), conclude. 
\end{proof}

\begin{remark}
Since, by Remark~\ref{cor:polynomial}, $q_t^{(1)}=O(t^{(s-a)/s})$, condition~\eqref{eq:bias-summability} follows from
$\sum_t t^{-a/s}\eta_t^+<+\infty$. Moreover, the simpler condition
$\sum_t\eta_t^-<+\infty$ 
implies~\eqref{eq:bias-square-summability}$(b)$, because
\(W_t\geq W_0\). In particular, for a bias $\eta_t\geq 0$ and such that 
\(\eta_t=O(t^{-\gamma})\), Assumption~\eqref{eq:bias-square-summability}$(b)$ trivially holds true; while  
conditions~\eqref{eq:bias-summability} and~\eqref{eq:bias-square-summability}$(a)$ hold whenever
$\gamma>\max\left\{1-a/s,1/2\right\}$. 
\end{remark}

\subsection{A randomly reinforced response-adaptive urn}
\label{subsubsec:response-adaptive:mart}

Consider a two-color randomly reinforced urn representing two treatments,
denoted by Black  ($B$) and White ($W$). Let $B_t$ and $W_t$ be their respective numbers of
balls,  with \(B_0,W_0>0\), and set $S_t:=B_t+W_t$. 
At time \(t+1\), a ball is sampled from the urn. Let \(I_{t+1}^W\) be the
indicator that color \(W\) is sampled, so that
\[
I_{t+1}^W\in\{0,1\},
\qquad
\Pbb(I_{t+1}^W=1\mid\F_t)=\frac{W_t}{S_t}.
\]
Moreover, let $C_{W,t}:=\sum_{s=1}^t I_s^W$, that is the number of allocations to treatment \(W\) up to time-step \(t\). 
The urn evolves according to
\begin{align}
W_{t+1}
&=W_t+R_{t+1}^W I_{t+1}^W,
\label{eq:RRU-W-update:mart}\\
B_{t+1}
&=B_t+R_{t+1}^B(1-I_{t+1}^W),
\end{align}
where \(R_{t+1}^W\) and \(R_{t+1}^B\) are, respectively, the nonnegative numbers of
white and black balls added when the corresponding color is drawn at time
\(t+1\).
Consequently,
\[
S_{t+1}
=
S_t+R_{t+1}^W I_{t+1}^W
+R_{t+1}^B(1-I_{t+1}^W).
\]
Assume that
\[
m_{B,t}:=\E(R_{t+1}^B\mid\F_t)\stackrel{a.s.}\longrightarrow m_B,
\qquad
m_{W,t}:=\E(R_{t+1}^W\mid\F_t)\stackrel{a.s.}\longrightarrow m_W\,,
\]
where \(0<m_W<m_B\) are deterministic. Thus \(W\) represents the asymptotically inferior treatment.
Assume that the non-negative reinforcements $[R_{t+1}^B,R_{t+1}^W]$ are, conditionally on the past $\F_t$, independent
of the draw $I_{t+1}^W$ and that 
\begin{equation}\label{eq:cond-rinfor}
E[(R^W_{t+1})^2\mid{\F}_t]=O(\,E[R^W_{t+1}\mid\F_t]\,)
\end{equation}
This last condition is  always satisfied when $(R_{t+1}^W)$ is bounded.

Randomly reinforced urns of this type arise in response-adaptive
clinical-trial designs \cite{BeCrPrRi10, CrimaldiLouisMinelli2022, CrimaldiLouisMinelli2023, Zhang2014Gaussian}. 

\begin{proposition} 
\label{prop:RRU-rate:mart} 
Assume that
\begin{equation}\label{eq:RRU-total-mass:mart}
\sum_t\frac1{t+1}
\left|
\frac{(t+1)m_{W,t}}{S_t}
-\frac{m_W}{m_B}
\right|
<+\infty
\qquad\text{a.s.}\,.
\end{equation}
Then there exists a strictly positive real random variable
$\widetilde{W}_\infty$ such that
\begin{equation}\label{eq:RRU-rate:mart}
t^{-m_W/m_B}W_t
\stackrel{a.s.}\longrightarrow
\widetilde{W}_\infty 
\end{equation}
and 
\begin{equation}\label{eq:RRU-allocation-rate:mart}
t^{-m_W/m_B}C_{W,t}
\stackrel{a.s.}\longrightarrow
\frac{\widetilde{W}_\infty}{m_W}.
\end{equation}
\end{proposition}

\begin{proof}
For \(t\geq0\), set $X_t:=W_t$, 
\[
p_{t+1}:=0,\qquad
r_{t+1}:=\frac1{t+1},\qquad
Y_{t+1}:=(t+1)R_{t+1}^WI_{t+1}^W,\qquad
\delta^*:=\frac{m_W}{m_B}.
\]
Then \(X^{(\pi)}_t=W_t\), \(\htr_{t+1}=1/(t+1)\), and
Assumption~\ref{ass:step-size} holds. Moreover,
\[
\E(Y_{t+1}\mid\F_t)
=(t+1)\frac{m_{W,t}W_t}{S_t}
=\delta^*W_t+\rho_t,
\qquad
\rho_t:=
W_t\left(\frac{(t+1)m_{W,t}}{S_t}-\delta^*\right).
\]
Thus Assumption~\ref{ass:reinforcement} holds. Conditional
independence and~\eqref{eq:cond-rinfor} give
\begin{multline*}
\Var(R_{t+1}^WI_{t+1}^W\mid\F_t)
\leq
\E\bigl((R_{t+1}^W)^2I_{t+1}^W\mid\F_t\bigr)
=
\frac{W_t}{S_t}
\E\bigl((R_{t+1}^W)^2\mid\F_t\bigr)
\\
=
O\!\left(
\E(R_{t+1}^WI_{t+1}^W\mid\F_t)
\right).
\end{multline*}
Since \(V_{t+1}=R_{t+1}^WI_{t+1}^W\), this verifies
Assumption~\ref{ass:moment}. 
Condition~\eqref{eq:RRU-total-mass:mart} says that
\[
\sum_t\frac{|\rho_t|}{(t+1)W_t}<+\infty\qquad\mbox{a.s.}\,.
\]
As \(W_0>0\) and \(W_t\) is non-decreasing, \(W_t>0\) for every \(t\).
Consequently,
\[
\sum_t\htr_{t+1}
 \left(\delta^*+\frac{\rho_t}{W_t}\right)
=
\delta^*\sum_t\frac1{t+1}
+\sum_t\frac{\rho_t}{(t+1)W_t}
=+\infty\quad\mbox{a.s.}\,,
\]
because the second series converges absolutely and \(\delta^*>0\).
Theorem~\ref{cor-X0}\textup{(a)} therefore gives
\(
W_t=X^{(\pi)}_t\stackrel{a.s.}\longrightarrow+\infty.
\) 
Also,
\[
\sum_t\htr_{t+1}\frac{|\rho_t|}{W_t}
=
\sum_t\frac{|\rho_t|}{(t+1)W_t}
<+\infty\qquad\mbox{a.s.}\,.
\]
The basic-recursion case of
Theorem~\ref{prop:absolute-relative-criterion}, followed by
Remark~\ref{cor:polynomial}, now proves
\eqref{eq:RRU-rate:mart}.\\
\indent It remains to consider the allocations. Since \(W_t\to+\infty\)
and \(W_t\) can change only when \(I_{t+1}^W=1\), necessarily
\(C_{W,t}\stackrel{a.s.}\to+\infty\). Put
\(
A_T^I:=\sum_{t=0}^{T-1}{W_t}/{S_t}.
\)
Apply Lemma~\ref{lem:compensated-increments} to the counting
process \(C_{W,T}\): its increments are \(I_{t+1}^W\) and
\[
\Var(I_{t+1}^W\mid\F_t)
\leq\E(I_{t+1}^W\mid\F_t)=\frac{W_t}{S_t}.
\]
Parts \ref{lem:compensated-increments:i} and \ref{lem:compensated-increments:ii} of the lemma, together with
\(C_{W,T}\stackrel{a.s.}\to+\infty\), show that
\[
A_T^I\stackrel{a.s.}\to+\infty,
\qquad
C_{W,T}\stackrel{a.s.}\sim A_T^I.
\]
Define the centered reinforcement martingale
\[
M_T^{W,\mathrm{cen}}
:=
\sum_{t=0}^{T-1}
I_{t+1}^W(R_{t+1}^W-m_{W,t}).
\]
Conditional independence, \eqref{eq:cond-rinfor}, and
\(m_{W,t}\to m_W\) imply
\[
\langle M^{W,\mathrm{cen}}\rangle_T
\leq
\sum_{t=0}^{T-1}
\frac{W_t}{S_t}
\E\bigl((R_{t+1}^W-m_{W,t})^2\mid\F_t\bigr)
=O(A_T^I).
\]
The martingale strong law therefore gives
\(M_T^{W,\mathrm{cen}}=o(A_T^I)=o(C_{W,T})\).
Since \(m_{W,t}\to m_W\) and \(C_{W,T}\stackrel{a.s.}\to+\infty\),
weighted Cesàro convergence gives
\[
\frac1{C_{W,T}}
\sum_{t=0}^{T-1}m_{W,t}I_{t+1}^W
\stackrel{a.s.}\longrightarrow m_W.
\]
Using
\(
W_T=W_0+
\sum_{t=0}^{T-1}m_{W,t}I_{t+1}^W
+M_T^{W,\mathrm{cen}},
\) 
we obtain \(W_T/C_{W,T}\stackrel{a.s.}\to m_W\). Combining this with
\eqref{eq:RRU-rate:mart} proves
\eqref{eq:RRU-allocation-rate:mart}.
\end{proof}

\begin{remark}
Condition \eqref{eq:RRU-total-mass:mart} required by Proposition~\ref{prop:RRU-rate:mart} is satisfied when  
$m_{B,t}=m_B$ and $m_{W,t}=m_W$ $\forall t\geq 0$, with $0<m_W<m_B$,
and $\sup_{t,\;c\in\{B,W\}} \E\bigl((R_{t+1}^c)^2\mid \F_t\bigr)\leq k$ a.s.\ for some $k>0$, see e.g.\ \cite[Theorem~3.2]{MayFlournoy2009} and \cite[Theorem~2.2]{ZhangHuCheungChan2014}.
 Indeed, in this setting, we can first prove that, for every \(\beta\in(m_W/m_B,1)\),
$W_t=o(t^\beta)$ and $ C_{W,t}=o(t^\beta)$ and,  consequently, condition \eqref{eq:RRU-total-mass:mart} holds.
\end{remark}

\subsection{Generalized Simon's model with signed score increments}
\label{subsec:signed-simon}\label{appl:simon}
Let \((\alpha_t)_t\) be a deterministic sequence such that $\alpha_0=1$, 
\(0<\alpha_t<1\) for \(t\geq1\), $\alpha_t\longrightarrow\alpha\in[0,1)$. Let \((S_t)_t\) be a (generalized) innovation process adapted to a filtration 
\((\F_t)_t\) and with ``innovation'' probabilities $(\alpha_t)$. At time \(t\), denote by \(\mathcal O_t\) the set of
colors already observed and set
$D_t:=\operatorname{card}(\mathcal O_t)$. 
At the next step, either a new color is observed or one of the colors
in \(\mathcal O_t\) is selected. More precisely,
\[
\Pbb(S_{t+1}=\text{``new''}\mid\F_t)=\alpha_t
\]
and the conditional probability of selecting an old color \(c\) at \(t\ge1\) is
\begin{equation}\label{eq:signed-simon-selection}
P_{t,c}
:=
\Pbb(S_{t+1}=c\mid\F_t)
=
(1-\alpha_t)
\frac{K_{t,c}+\beta}{t+\beta D_t}\,,
\end{equation}
where $K_{t,c}$ is a non-negative \(\F_t\)-measurable score assigned to color $c$ at time-step $t$ such that 
\begin{equation}\label{eq:signed-simon-normalization}
K_{t,c}+\beta>0\qquad\mbox{and}\qquad 
\sum_{c\in\mathcal O_t}K_{t,c}=t.
\end{equation}
These conditions assure that $0<P_{t,c}\leq 1$ and 
$\sum_{c\in\mathcal O_t}P_{t,c}=1-\alpha_t$, so that \eqref{eq:signed-simon-selection}, together with the innovation
probability \(\alpha_t\), defines a probability distribution. When a new color appears, it is added to \(\mathcal O_t\) with initial
score equal to one (and so, by~\eqref{eq:signed-simon-normalization}, $\beta>-1$), while the scores of the old colors remain unchanged. When
an old color is selected, the scores may undergo signed updates: that is, for
\(c\in\mathcal O_t\), we set 
\begin{equation}\label{eq:signed-simon-score-update}
K_{t+1,c}:=K_{t,c}+Z_{t+1,c}\,,
\end{equation}
where the update \(Z_{t+1,c}\) may be negative. Since $K_{t+1,c}\geq 0$ and \eqref{eq:signed-simon-normalization}, for every
\(c\in\mathcal O_t\), we require 
\[
K_{t,c}+Z_{t+1,c}\geq0,
\qquad
K_{t,c}+Z_{t+1,c}+\beta>0,
\qquad
\sum_{c\in\mathcal O_t}Z_{t+1,c}
=\mathbf1_{\{S_{t+1}\in\mathcal O_t\}}.
\]
Moreover, 
the signed updates are assumed to preserve the usual Simon drift:
\begin{equation}\label{eq:signed-simon-conditional-mean}
\E(Z_{t+1,c}\mid\F_t)
=
P_{t,c}
=
(1-\alpha_t)
\frac{K_{t,c}+\beta}{t+\beta D_t}>0.
\end{equation}
Thus an individual score may decrease, although its conditional mean
increment is strictly positive. Assume also that
\begin{equation}\label{eq:signed-simon-variance}
\Var(Z_{t+1,c}\mid\F_t)
=O\bigl(\E(Z_{t+1,c}\mid\F_t)\bigr).
\end{equation}
For a fixed color \(c\), define its first-appearance time as 
\[
\tau_c:=\inf\{t\geq1:c\in\mathcal O_t\},
\]
which is a stopping time since
\(\{\tau_c\leq t\}=\{c\in\mathcal O_t\}\in\F_t\).
Note that, when $-1<\beta\leq 0$, for each observed color $c$, that is for each color $c$ with $\tau_c<+\infty$, we have 
$K_{t,c}>0$ for all $t\geq \tau_c$. Indeed, by the first condition in~\eqref{eq:signed-simon-normalization}, we have 
$K_{t,c}>-\beta\geq 0$. 
For the following result, 
when \(-1<\beta<0\), we need to assume 
\begin{equation}\label{eq:signed-simon-beta-nondegeneracy}
\liminf_{t}
\left(1+\frac{\beta}{K_{t,c}}\right)>0
\qquad\text{a.s. on }\{\tau_c<+\infty\},
\end{equation}
while, when \(\beta=0\), we need to assume 
\begin{equation}\label{eq:signed-simon-loss}
\sum_{t\geq\tau_c}\frac{Z_{t+1,c}^-}{K_{t,c}}
<+\infty
\qquad\text{a.s. on }\{\tau_c<+\infty\}.
\end{equation}
The generalized Simon model in \cite{crimaldi2026triggeredurnmodelsfrequently} is recovered by taking 
\(
Z_{t+1,c}
=
\ind_{\{S_{t+1}=c\}}, 
\) 
so that \(K_{t,c}\) corresponds to the number of observations of color \(c\) (in the classical Simon model \cite{simon-standard} we have in addition $\alpha_t=\alpha$ for all $t\geq 1$ and $\beta=0$). 
In this case 
\eqref{eq:signed-simon-variance} and
\eqref{eq:signed-simon-loss} are automatic. On the other hand, when an old color is selected, an independent fair
signed transfer of size
\(2^{-(t+2)}\min_{i\in\{1,2\}}\{K_{t,c_i},K_{t,c_i}+\beta\}\)
between the first two observed colors gives an admissible update
with genuinely negative increments.
 
\begin{proposition}
\label{prop:signed-innovation-fixed-color}
Under the preceding assumptions,   
on $\{\tau_c<+\infty\}$,
\begin{align}
\frac{\ln K_{t,c}}{\ln t}
&\stackrel{a.s.}\longrightarrow
\lambda:=\frac{1-\alpha}{1+\beta\alpha}>0,
\label{eq:signed-innovation-first-order-K}\\
\frac{\ln P_{t,c}}{\ln t}
&\stackrel{a.s.}\longrightarrow
\lambda-1
=-\frac{\alpha(1+\beta)}{1+\beta\alpha}.
\label{eq:signed-innovation-first-order-P}
\end{align}
If, in addition,
\begin{equation}\label{eq:signed-innovation-alpha-summability}
\sum_t\frac{|\alpha_t-\alpha|}{t+1}<+\infty,
\end{equation}
then there is a strictly positive real random variable
\(\widetilde{K}_{\infty,c}\) such that,
on $\{\tau_c<+\infty\}$,
\begin{align}
t^{-\lambda}K_{t,c}
&\stackrel{a.s.}\longrightarrow \widetilde{K}_{\infty,c},
\label{eq:signed-innovation-exact-K}\\
t^{1-\lambda}P_{t,c}
&\stackrel{a.s.}\longrightarrow
\lambda \widetilde{K}_{\infty,c}.
\label{eq:signed-innovation-exact-P}
\end{align}
\end{proposition}

\begin{proof}
Set $\kappa:=1+\beta\alpha>0$ and $\delta^*:=1-\alpha>0$ so that $\lambda={\delta^*}/{\kappa}$ and let 
$r_{t+1}:=1/{(\kappa t+1)}$. Moreover, set 
\[
J_s:=\mathbf1_{\{S_s=\text{``new''}\}},
\qquad
M_t^D:=\sum_{s=1}^t(J_s-\alpha_{s-1}).
\]
Since \(D_t=\sum_{s=1}^tJ_s\), the bounded-increment
martingale strong law gives
\begin{equation}\label{eq:signed-innovation-D-LLN}
\frac{D_t}{t}\stackrel{a.s.}\longrightarrow\alpha.
\end{equation}
Define 
\[
W_{t,c}:=K_{t,c}+\beta,\qquad
A_t:=(1-\alpha_t)
       \frac{\kappa t+1}{t+\beta D_t}.
\]
Then, on $\{\tau_c<+\infty\}$, we have  \(W_{t,c}>0\) after the first appearance of color $c$, i.e. for $t\geq \tau_c$, 
\(A_t\stackrel{a.s.}\to\delta^*\), and
\begin{equation}\label{eq:signed-simon-rho}
P_{t,c}=r_{t+1}A_tW_{t,c},
\qquad
A_tW_{t,c}
=\delta^*W_{t,c}+\rho^W_{t,c},
\qquad
\rho^W_{t,c}:=(A_t-\delta^*)W_{t,c}.
\end{equation}
We first check that \(W_{t,c}\) stays bounded away
from zero on \(\{\tau_c<+\infty\}\). This is immediate
when \(\beta>0\), since \(W_{t,c}\geq\beta\).
When \(\beta=0\), the model gives
\(0\leq U_{t+1,c}:=Z_{t+1,c}^-/K_{t,c}<1\).
Condition~\eqref{eq:signed-simon-loss} makes
\(\prod_{t\geq\tau_c}(1-U_{t+1,c})\) strictly
positive; hence
\[
K_{T,c}\geq K_{\tau_c,c}
\prod_{t=\tau_c}^{T-1}(1-U_{t+1,c})\quad\mbox{with } K_{\tau_c,c}=1
\]
is bounded away from zero. Finally, when
\(-1<\beta<0\), we have \(K_{t,c}>-\beta\), and
\eqref{eq:signed-simon-beta-nondegeneracy} implies
\(W_{t,c}=K_{t,c}(1+\beta/K_{t,c})\) is bounded
away from zero eventually. The finitely many
earlier values are strictly positive so that 
\begin{equation}\label{eq:simon-shift-lower-bound}
\inf_{t\geq\tau_c}W_{t,c}>0
\qquad\text{a.s. on }\{\tau_c<+\infty\}.
\end{equation}
For clarity about the application of the main
theorems, we use the following localization.
Fix deterministic \(n,m\geq1\), set
\(G_n:=\{\tau_c\leq n\}\in\F_n\), and, on \(G_n\),
write \(\Sigma_t:=\Var(Z_{t+1,c}\mid\F_t)\) and
\[
\zeta_{n,m}:=
\inf\{t\geq n:
 W_{t,c}<m^{-1}\ \text{or}\ \Sigma_t>mP_{t,c}\}.
\]
Stop \(W_{t,c}\) at \(\zeta_{n,m}\), and set the
stopped process equal to zero off \(G_n\).
Before stopping, take
\[
Y_{t+1}:=A_tW_{t,c},
\qquad
\Delta N_{t+1}:=
\frac{Z_{t+1,c}-P_{t,c}}
     {r_{t+1}W_{t,c}};
\]
after stopping, take both quantities to be zero.
Since \(Y_{t+1}=A_tW_{t,c}\) is \(\F_t\)-measurable,
\(\Delta M^Y_{t+1}=0\) on active steps. Thus the martingale term in the recursion is
\(r_{t+1}W_{t,c}\Delta N_{t+1}=Z_{t+1,c}-P_{t,c}\).
This gives the modified recursion
\eqref{eq:recursion:mart}, in shifted time, with
\(p_{t+1}=0\). On active steps its increment has
conditional mean \(P_{t,c}\) and conditional
variance at most \(mP_{t,c}\).
On active steps, \(W_{t,c}\geq m^{-1}\),
\(\Sigma_t\leq mP_{t,c}\), and \(P_{t,c}\leq1\).
Hence
\[
\E\bigl((\Delta N_{t+1})^2\mid\F_t\bigr)
=\frac{\Sigma_t}{r_{t+1}^2W_{t,c}^2}
\leq\frac{m^3}{r_{t+1}^2}.
\]
For fixed \(n,m\), this deterministic bound is finite at
each time, so the localized increments define the
martingale \((N_t)_t\) required in~\eqref{eq:recursion:mart}.
Summability of its relative increments is verified below.
By \eqref{eq:signed-simon-variance} and \eqref{eq:simon-shift-lower-bound}, almost every
path in \(\{\tau_c<+\infty\}\) belongs to
\(G_n\cap\{\zeta_{n,m}=+\infty\}\) for some
deterministic \(n,m\). On each such event the
stopped process is \(W_{t,c}\), and the effective
step sizes satisfy Assumption~\ref{ass:step-size}.
Also,
\[
\frac{(\rho^W_{t,c})^-}{W_{t,c}}
=(\delta^*-A_t)^+\stackrel{a.s.}\longrightarrow0<\delta^*.
\]
Theorem~\ref{cor-X0}\ref{cor-X0:b} therefore gives
the alternatives \(W_{t,c}\to0\) and
\(W_{t,c}\to+\infty\). The former is excluded by
\(W_{t,c}\geq m^{-1}\). Taking the countable union
over \(n,m\), we obtain
\begin{equation}\label{eq:signed-simon-score-diverges}
W_{t,c}\stackrel{a.s.}\longrightarrow+\infty,
\qquad
K_{t,c}\stackrel{a.s.}\longrightarrow+\infty
\quad\text{on }\{\tau_c<+\infty\}.
\end{equation}
On each event just considered, by~\eqref{eq:signed-simon-rho} we have
\(P_{t,c}/W_{t,c}=r_{t+1}A_t\stackrel{a.s.}\to 0\).
Lemma~\ref{lem:compensated-increments}
\ref{lem:compensated-increments:iii}, applied to
the stopped process, yields
\[
\sum_t
\frac{\Var(Z_{t+1,c}\mid\F_t)}
     {W_{t,c}^2}<+\infty.
\]
Thus the relative martingale increment
\[
\Delta\mathcal M_{t+1}
=
\frac{Z_{t+1,c}-P_{t,c}}
     {(1+\delta^*r_{t+1})W_{t,c}}
\]
satisfies Assumption~\ref{ass:relative-noise:mart2}.
Furthermore, \(A_t\stackrel{a.s.}\to\delta^*\) and
\(\sum_t r_{t+1}^2<+\infty\), so
\[
\sum_t\left(
r_{t+1}\frac{(\rho^W_{t,c})^+}{W_{t,c}}
\right)^2
=
\sum_t\bigl(r_{t+1}(A_t-\delta^*)^+\bigr)^2
<+\infty.
\]
This verifies
Assumption~\ref{ass:quadratic-relative-perturbation}.
On each event \(G_n\cap\{\zeta_{n,m}=+\infty\}\),
the perturbation sum in~\eqref{eq:Gamma}, for the
localized recursion started at \(n\), is
\[
\Gamma_{T-1}
=\sum_{t=n}^{T-1}r_{t+1}
       \frac{\rho^W_{t,c}}{W_{t,c}}
=\sum_{t=n}^{T-1}r_{t+1}(A_t-\delta^*).
\]
Theorem~\ref{thm:main-positive-all-maringale-noise}
therefore gives
\[
\ln\bigl(q_T^{(\delta^*)}W_{T,c}\bigr)
-\Gamma_{T-1}
\stackrel{a.s.}\longrightarrow L
\]
for a real random variable \(L\) on
\(\{\tau_c<+\infty\}\). Since
\(A_t-\delta^*\stackrel{a.s.}\to0\) and
\(\sum_{t<T}r_{t+1}
=\kappa^{-1}\ln T+O(1)\), we have
\(\Gamma_{T-1}=o(\ln T)\).
Also,
\(\ln q_T^{(\delta^*)}
=-\lambda\ln T+O(1)\).
Hence
\(\ln W_{T,c}/\ln T\stackrel{a.s.}\to\lambda\).
Because \(W_{T,c}=K_{T,c}+\beta\stackrel{a.s.}\to+\infty\),
this proves
\eqref{eq:signed-innovation-first-order-K}.
The identity \(P_{t,c}=r_{t+1}A_tW_{t,c}\)
presented in~\eqref{eq:signed-simon-rho} then gives
\eqref{eq:signed-innovation-first-order-P}.
\\
\indent Finally, assume
\eqref{eq:signed-innovation-alpha-summability}.
The bound \(\E((M_t^D)^2)\leq t/4\) and Cauchy–Schwarz give 
$\sum_t \E[\,|M_t^D|\,]/t^2=O(\sum_t t^{-3/2})<+\infty$. 
Moreover, an exchange of the order of
summation in
\(D_t-\alpha t
=M_t^D+\sum_{s<t}(\alpha_s-\alpha)\)
gives by~\eqref{eq:signed-innovation-alpha-summability}
\begin{equation}\label{eq:signed-simon-D-summability}
\sum_{t\geq1}\frac{|D_t-\alpha t|}{t^2}
<+\infty\qquad\text{a.s.}
\end{equation}
Then, since 
\[
|A_t-\delta^*|
=O\left(
|\alpha_t-\alpha|
+\frac{|D_t-\alpha t|+1}{t}
\right)\,,
\] we get 
\begin{equation}\label{eq:signed-simon-a-summability}
\sum_{t\geq1}\frac{|A_t-\delta^*|}{t}
<+\infty\qquad\text{a.s.}
\end{equation}
Consequently, 
\[
\sum_t r_{t+1}\frac{|\rho^W_{t,c}|}{W_{t,c}}
=
\sum_t r_{t+1}|A_t-\delta^*|
<+\infty
\quad\text{on }\{\tau_c<+\infty\}.
\]
Theorem~\ref{prop:absolute-relative-criterion}
gives a strictly positive real limit for
\(q_T^{(\delta^*)}W_{T,c}\) on this event.
Remark~\ref{cor:polynomial}, applied with
\(p_t=0\) and
\(r_t=\kappa^{-1}t^{-1}+O(t^{-2})\), gives
\(q_T^{(\delta^*)}/T^{-\lambda}\to\ell\)
for a deterministic \(\ell\in(0,+\infty)\).
Thus, on $\{\tau_c<+\infty\}$, \(T^{-\lambda}W_{T,c}\) has a strictly
positive real limit
\(\widetilde K_{\infty,c}\).
Since \(\lambda>0\) and \(K_{T,c}=W_{T,c}-\beta\),
this is \eqref{eq:signed-innovation-exact-K}.
Finally, on $\{\tau_c<+\infty\}$
\[
T^{1-\lambda}P_{T,c}
=(Tr_{T+1})A_T(T^{-\lambda}W_{T,c})
\stackrel{a.s.}\longrightarrow
\frac{\delta^*}{\kappa}\widetilde K_{\infty,c}
=\lambda\widetilde K_{\infty,c},
\]
proving \eqref{eq:signed-innovation-exact-P}.
\end{proof}

\subsection{Branching Process in a Varying Environment (BPVE)}
\label{subsec:applications-signed-increments}
 The following example deals with a 
model which is a member of the class of branching processes in a varying
environment studied in \cite{Kersting2020}.
Fix $0<\eta<a<1$, $b>0$ 
and choose \(c>0\) so large that the probabilities \(\alpha_t+\beta_t\leq1\) and \(a-bc^{-\eta}>0\).  Let $\{(\xi_{t,j})_t: j\geq 1\}$ be a collection of stochastic processes, adapted to a filtration $\F$, and such that, 
conditionally on \(\F_t\), the variables
\(\{\xi_{t+1,j}:j\geq 1\}\) are independent and
\begin{equation}
\label{eq:bpve-rare-death-offspring}
 \Pbb(\xi_{t+1,j}=2\mid\F_t)=\alpha_t:=\frac{a}{t+c},\qquad
 \Pbb(\xi_{t+1,j}=0\mid\F_t)=\beta_t:=\frac{b}{(t+c)^{1+\eta}},
\end{equation}
with the remaining probability $(1-\alpha_t-\beta_t)$ assigned to \(\xi_{t+1,j}=1\).
Define the non-negative stochastic process \((Z_t)_t\), adapted to the filtration~$\F$, setting 
\[
Z_0\in\mathbb{N} >0,\qquad Z_{t+1}:=\sum_{j=1}^{Z_t}\xi_{t+1,j}\,.
\]
This is a BPVE in the terminology of \cite{Kersting2020}. We define its survival event as  
\[
\mathcal S:=\liminf_t\{Z_t>0\}.
\]
 Notice that \(\Pbb(\mathcal S)>0\): indeed,
the probability that one fixed ancestral line never produces zero offspring
is bounded below by the product 
\(
 \prod_t(1-b(t+c)^{-1-\eta})>0
\).  The following result shows that, for the considered  BPVE, on the survival event $\mathcal S$ the size of the population exhibits a power-law growth, although
 the  increment of the size between two consecutive generation can be strictly negative infinitely often. The strict positive normalized limit also follows from
\cite[Theorem~2]{Kersting2020}; here we additionally show that
population decreases occur infinitely often on the survival event.
\begin{proposition}
\label{prop:bpve-signed-application}
On $\mathcal S$, we have 
\begin{equation}\label{eq:thesiZ}
 t^{-a}\, Z_t\stackrel{a.s.}\longrightarrow \widetilde{Z}_\infty\,,
\end{equation}
where $\widetilde{Z}_\infty$ is a strictly positive real random variable. Moreover, 
the effective increment is genuinely signed: that is,
\begin{equation}
\label{eq:thesisZ-2}
 \Pbb\bigl( Z_{t+1}-Z_t<0\ \text{infinitely often}
          \mid\mathcal S\bigr)=1\,.
\end{equation}
\end{proposition}

\begin{proof}
Set $\widetilde V_{t+1}:=Z_{t+1}-Z_t$. 
Then
\[
\mu_t
:=
\E(\widetilde V_{t+1}\mid\F_t)
=
Z_t(\alpha_t-\beta_t)
=
\frac1{t+c}
\left(
 aZ_t-\frac{bZ_t}{(t+c)^\eta}
\right),
\]
and
\[
\Var(\widetilde V_{t+1}\mid\F_t)
\leq
Z_t(\alpha_t+\beta_t)
=O(\mu_t),
\]
because \(\sup_t\beta_t/\alpha_t<1\). Setting
\[
p_{t+1}:=0,
\qquad
r_{t+1}:=\frac1{t+c},
\qquad
\delta^*:=a,
\qquad
\rho_t:=-\frac{bZ_t}{(t+c)^\eta},
\]
$$
Y_{t+1}:=(t+c)\mu_t=\delta^*Z_t+\rho_t\qquad\mbox{and}\qquad 
\Delta N_{t+1}
:=
\frac{(t+c)(\widetilde V_{t+1}-\mu_t)}{Z_t}
\ind_{\{Z_t>0\}},
$$
we obtain the recursion~\eqref{eq:recursion:mart} for $(Z_t)_t$.  
The step-size
assumptions~\ref{ass:step-size} are immediate.
On the survival event \(\mathcal S\), the process \((Z_t)\) takes values in $\mathbb{N}\setminus\{0\}$ eventually, and
\[
\frac{\rho_t^-}{Z_t}
=\frac{b}{(t+c)^\eta}\longrightarrow0<\delta^*.
\]
The conditional-variance bound above and
Theorem~\ref{cor-X0}\ref{cor-X0:b} show that, on
\(\mathcal S\), either \(Z_t\to0\) or \(Z_t\to+\infty\).
The first alternative is impossible since \(Z_t\geq 1\) eventually 
on \(\mathcal S\). Hence
\(
Z_t\stackrel{a.s.}\longrightarrow+\infty\) on $\mathcal S$.
Since ${\mu_t}/{Z_t}=\alpha_t-\beta_t\longrightarrow 0$, 
Lemma~\ref{lem:compensated-increments}\ref{lem:compensated-increments:iii} gives
\[
\sum_t
\frac{\Var(\widetilde V_{t+1}\mid\F_t)}{Z_t^2}
<+\infty
\qquad\text{a.s. on }\mathcal S.
\]
As \(Y_{t+1}\) is $\F_t$-measurable, on \(\mathcal S\) we have 
\[
\Delta\mathcal M_{t+1}
=
\frac{\widetilde V_{t+1}-\mu_t}
 {Z_t(1+a/(t+c))}
\]
and so Assumption~\ref{ass:relative-noise:mart2} holds. Furthermore,
\[
\sum_t
\htr_{t+1}\frac{|\rho_t|}{Z_t}
=
b\sum_t\frac1{(t+c)^{1+\eta}}
<+\infty.
\]
Theorem~\ref{prop:absolute-relative-criterion} and
Remark~\ref{cor:polynomial} therefore give~\eqref{eq:thesiZ}.
\\
\indent It remains to prove that negative increments occur infinitely often.
Let \(A_{t+1}\) be the event that exactly one individual has no
offspring and all the others have one offspring. Then
\(
A_{t+1}\subseteq\{\widetilde V_{t+1}=-1\}
\)
and
\[
\Pbb(A_{t+1}\mid\F_t)
=
Z_t\beta_t(1-\alpha_t-\beta_t)^{Z_t-1}.
\]
By \eqref{eq:thesiZ}, $Z_t\beta_t\stackrel{a.s.}\sim b\widetilde Z_\infty t^{a-1-\eta}$ on $\mathcal S$.
Moreover,
\(
Z_t(\alpha_t+\beta_t)=O(t^{a-1})\to 0,
\)
so
\(
(1-\alpha_t-\beta_t)^{Z_t-1}\stackrel{a.s.}\to 1
\).  Thus, we have 
\(
\Pbb(A_{t+1}\mid\F_t)
\stackrel{a.s.}\sim
b\widetilde Z_\infty t^{a-1-\eta}
\) on $\mathcal S$. 
Since \(a>\eta\), the conditional probabilities have divergent sum.
Lévy's extension of Borel--Cantelli lemma \cite{will91} proves~\eqref{eq:thesisZ-2}.
\end{proof}

\subsection{Non-negative linear stochastic approximation with
multiplicative noise}
\label{subsubsec:linear-SA:mart}

Let \((\gamma_t)_t\) be a deterministic non-negative sequence satisfying
\begin{equation}\label{eq:linear-SA-general-step}
\gamma_t\longrightarrow0,
\qquad
\sum_t\gamma_t=+\infty,
\qquad
\sum_t\gamma_t^2<+\infty.
\end{equation}
Consider
\begin{equation}\label{eq:linear-SA-update:mart}
X_0>0,\qquad X_{t+1}
=
\left(
1-\gamma_{t+1}
+\sigma\gamma_{t+1}\xi_{t+1}
\right)X_t,
\end{equation}
where \(\sigma\geq 0\) and \((\xi_{t+1})_t\) is a
square-integrable martingale-difference sequence such that
\begin{equation}\label{eq:linear-SA-noise-L2}
\E(\xi_{t+1}\mid\F_t)=0,
\qquad
\sup_t\E(\xi_{t+1}^2\mid\F_t)<+\infty
\qquad\text{a.s.}
\end{equation}
Assume that
\begin{equation}\label{eq:linear-SA-positivity}
1-\gamma_{t+1}
+\sigma\gamma_{t+1}\xi_{t+1}>0
\qquad\text{a.s. for every }t.
\end{equation}
Hence $X_t>0$ a.s. for every $t$. 
Note that since the quantity in \eqref{eq:linear-SA-positivity}
is strictly positive and has conditional expectation
\(1-\gamma_{t+1}\), we have \(\gamma_{t+1}<1\) for
every \(t\). Thus the full product \(q_t^{(1)}\)
appearing below is well defined.

The following result complements the concentration theory of
\cite{ChenMaguluriZubeldia2025} by identifying the exact a.s.\
normalization and the strict positivity of its random coefficient.

\begin{proposition}
\label{prop:linear-SA:mart} There exists a strictly positive real random variable $X^*_\infty$ such that 
\begin{equation}\label{eq:thesis-SA}
q_t^{(1)}X_t
\stackrel{a.s.}\longrightarrow
X_\infty^*,\qquad q_t^{(1)}
:=
\prod_{s=0}^{t-1}(1-\gamma_{s+1})^{-1}, 
\end{equation}
and, consequently, there exists a strictly positive real random variable $\widetilde{X}_\infty$ such that 
\begin{equation}\label{eq:linear-SA-general-rate}
e^{{g_t}}X_t
\stackrel{a.s.}\longrightarrow
\widetilde{X}_\infty, \qquad\mbox{with}\qquad {g_t}:=\sum_{s=0}^{t-1}\gamma_{s+1}.
\end{equation}
\end{proposition}
\begin{proof}
Choose \(t_0\) such that
\(4\gamma_{t+1}\leq1\) for \(t\geq t_0\). For $t\geq t_0$, set
\[
p_{t+1}:=2\gamma_{t+1},\quad
r_{t+1}:=\gamma_{t+1},\quad
Y_{t+1}:=X_t,\quad
\delta^*:=1,\quad
\rho_t:=0,\quad
\Delta N_{t+1}:=\sigma\xi_{t+1}.
\]
Then \eqref{eq:linear-SA-update:mart} has the
form~\eqref{eq:recursion:mart}, and
Assumption~\ref{ass:reinforcement} holds. Moreover,
\[
\htr_{t+1}
=\frac{\gamma_{t+1}}{1-2\gamma_{t+1}}\,,
\]
so \eqref{eq:linear-SA-general-step} verifies
Assumption~\ref{ass:step-size}. Moreover, define
\[
\pi_{t_0}:=1,\qquad
\pi_T:=
\prod_{t=t_0}^{T-1}(1-2\gamma_{t+1})^{-1}
\quad(T>t_0),
\qquad
X_T^{(\pi)}:=\pi_TX_T.
\]
We have 
\(
\widetilde V_{t+1}
=X_{t+1}^{(\pi)}-X_t^{(\pi)}
=\htr_{t+1}X_t^{(\pi)}
 (1+\sigma\xi_{t+1}).
\)
Consequently,
\[
\frac{\Var(\widetilde V_{t+1}\mid\F_t)}
     {(X_t^{(\pi)})^2}
=\sigma^2(\htr_{t+1})^{2}
  \E(\xi_{t+1}^2\mid\F_t).
\]
Since \(\htr_{t+1}\leq2\gamma_{t+1}\) for $t\geq t_0$, \eqref{eq:linear-SA-general-step}
and~\eqref{eq:linear-SA-noise-L2} give
\[
\sum_{t\geq t_0}
\frac{\Var(\widetilde V_{t+1}\mid\F_t)}
     {(X_t^{(\pi)})^2}
\leq
4\sigma^2
\left(\sup_t\E(\xi_{t+1}^2\mid\F_t)\right)
\sum_{t\geq t_0}\gamma_{t+1}^2
<+\infty
\qquad\text{a.s.}
\]
Because $X_t>0$ for every $t$ and \(\rho_t=0\),
the event \(E_c\) of Theorem~\ref{cor-X0}\ref{cor-X0:c} has probability one. That theorem
therefore gives
\(
X_T^{(\pi)}\stackrel{a.s.}\longrightarrow+\infty.
\) 
It remains to check the hypotheses
of Theorem~\ref{prop:absolute-relative-criterion}.
Since \(Y_{t+1}=X_t\) is predictable,
\[
\Delta\mathcal M_{t+1}
=\frac{\htr_{t+1}}{1+\htr_{t+1}}
  \Delta N_{t+1}
=\frac{\sigma\gamma_{t+1}\xi_{t+1}}
       {1-\gamma_{t+1}}.
\]
The denominator is bounded away from zero for $t\geq t_0$. Thus \eqref{eq:linear-SA-general-step}
and~\eqref{eq:linear-SA-noise-L2} imply
\[
\sum_{t\geq t_0}
\E\bigl((\Delta\mathcal M_{t+1})^2
        \mid\F_t\bigr)
<+\infty
\qquad\text{a.s.},
\]
which verifies Assumption~\ref{ass:relative-noise:mart2}.
Also,
\(
\sum_{t\geq t_0}
\htr_{t+1}{|\rho_t|}/{X_t}=0.
\)
Theorem~\ref{prop:absolute-relative-criterion} can be applied 
 and yields the a.s. convergence of the tail-normalized process
\( (q_T^{(1)}/q_{t_0}^{(1)})X_T\) toward a strictly positive real 
 random variable and, 
multiplying by the positive deterministic factor
\(q_{t_0}^{(1)}\), we get \eqref{eq:thesis-SA}. 
Finally, \(\gamma_t<1\) and \(\sum_t\gamma_t^2<+\infty\) 
 imply
\[
\ln q_T^{(1)}=g_T+d+o(1),
\qquad
d:=
\sum_{s\geq0}
\bigl[-\ln(1-\gamma_{s+1})-\gamma_{s+1}\bigr]
<+\infty.
\]
Hence \(q_T^{(1)}e^{-g_T}\to e^d\), and
\eqref{eq:linear-SA-general-rate} follows with
\(\widetilde X_\infty=e^{-d}X_\infty^*>0\).
\end{proof}

\noindent {\bf Author contributions}\\
All the authors contributed equally to the present work.

\appendix

\section{Existence of a finite  limit for the rescaled process}\label{sec:upper}

 The following result is an immediate consequence of Robbins--Siegmund Theorem \cite{RobbinsSiegmund1971}.

\begin{theorem}[Finite normalized limit]\label{thm:finite-limit}
Consider either the basic dynamics \eqref{eq:recursion} or the
modified dynamics \eqref{eq:recursion:mart}, together with
Assumptions~\ref{ass:reinforcement} and~\ref{ass:lim_sup2}.
Then, the rescaled process $X^*_t=q_t^{(\delta^*)}X_t$ converges almost surely to a non-negative real random variable
$X^*_\infty$ and $\sum_{t}q_{t+1}^{(\delta^*)}r_{t+1}|\rho_t|< +\infty$ a.s.\,.
 Consequently, we can write $X_t=O\bigl(1/q_t^{(\delta^*)}\bigr)$. 
\end{theorem}

\begin{proof}
For both recursions, we have 
\[
E[X_{t+1}^*\mid{\F}_t]
=X_t^* +q_{t+1}^{(\delta^*)}r_{t+1}\rho_t=X_t^* +q_{t+1}^{(\delta^*)}r_{t+1}\rho^+_t - q_{t+1}^{(\delta^*)}r_{t+1}\rho^-_t .
\]
Hence, $(X^*_t)$
 is a non-negative almost super-martingale.  Thus, the Robbins--Siegmund Theorem \cite{RobbinsSiegmund1971} can be applied 
  and the thesis is proven. 
\end{proof}

\section{Technical lemma}\label{sec:technical}
The following lemma extends  
\cite[Theorem~46, p.~40]{del-mey} to stochastic processes with possibly unbounded,
signed increments, under a conditional-variance bound in terms of their
non-negative conditional means.
\begin{lemma}
\label{lem:compensated-increments}
Let \((Z_t)_t\) be a real-valued process adapted to
\((\F_t)_t\), whose increments have finite conditional second
moments almost surely. 
Suppose that 
\begin{equation}\label{eq:lemma-cond-var-bound}
\Var(\Delta Z_{t+1}\mid\F_t)
=O\bigl(\E(\Delta Z_{t+1}\mid\F_t)\bigr).
\end{equation}
Define
\[
A_T:=\sum_{t\leq T-1}\E(\Delta Z_{t+1}\mid\F_t),
\qquad
M_T^Z:=\sum_{t\leq T-1} (\Delta Z_{t+1}-\E(\Delta Z_{t+1}\mid\F_t)).
\]
Then:
\begin{enumerate}[label=\textup{(\roman*)}]
\item\label{lem:compensated-increments:i} On \(\{A_\infty<+\infty\}\), both \(M_T^Z\) and \(Z_T\)
converge a.s.\ to a real random variable.
\item\label{lem:compensated-increments:ii} On \(\{A_T\to+\infty\}\), $M_T^Z=o(A_T)$ and $Z_T\stackrel{a.s.}\sim A_T$.
\item\label{lem:compensated-increments:iii} On 
$$
\{Z_t\to+\infty\}
\cap
\left\{\frac{\E(\Delta Z_{t+1}\mid\F_t)}{Z_t}\to0\right\},
$$
we have~\footnote{The sums are taken starting from a sufficiently large random time
(inclusive) after which \(Z_t>0\).}
\[
\sum_t\frac{\E(\Delta Z_{t+1}\mid\F_t)}{Z_t^2}<+\infty\qquad\mbox{and}\qquad 
\sum_t
\frac{\Var(\Delta Z_{t+1}\mid\F_t)}{Z_t^2}
<+\infty
\quad\text{a.s.}
\]
\end{enumerate}
\end{lemma}

\begin{proof}
The centered process \(M^Z\) is a locally square-integrable
martingale, and
\[
Z_T=Z_0+A_T+M_T^Z,\quad
\langle M^Z\rangle_T
=\sum_{t<T}\E((\Delta M^Z_{t+1})^2\mid\F_t) 
=\sum_{t<T}\Var(\Delta Z_{t+1}\mid\F_t).
\]
By the hypotheses, almost everywhere, there are an integer-valued random time \(\tau\) 
(not necessarily a stopping time) and a non-negative real random variable  \(C\)
such that, for \(t\geq\tau\),
\[
0\leq\E(\Delta Z_{t+1}\mid\F_t),
\qquad
\Var(\Delta Z_{t+1}\mid\F_t)
\leq C\,\E(\Delta Z_{t+1}\mid\F_t)
\]
and, for
\(T\geq\tau+1\),
\[
\begin{aligned}
\langle M^Z\rangle_T-\langle M^Z\rangle_\tau
=\sum_{t=\tau}^{T-1}
  \Var(\Delta Z_{t+1}\mid\F_t)\leq C\sum_{t=\tau}^{T-1}
  \E(\Delta Z_{t+1}\mid\F_t)
 =C(A_T-A_\tau).
\end{aligned}
\]
Thus \(A_T\) is eventually non-decreasing and
\(A_\infty:=\lim_{T\to\infty}A_T\) exists almost surely in
\(\mathbb R\cup\{+\infty\}\). 
Consequently, on \(\{A_\infty<+\infty\}\) we have
\(\langle M^Z\rangle_\infty<+\infty\); while for the same relation on
\(\{A_T\to+\infty\}\) we have by \eqref{eq:lemma-cond-var-bound}  
\(\langle M^Z\rangle_T=O(A_T)\). Then, 
on \(\{A_\infty<+\infty\}\), 
Localized Martingale Convergence Theorem \cite[Theorem~2.3.4, p.~46]{lapeyre2024montecarlo}
gives
convergence of \(M_T^Z\), and hence of \(Z_T\), to a real random variable. This proves
\ref{lem:compensated-increments:i}.
On \(\{A_T\to+\infty\}\), if
\(\langle M^Z\rangle_\infty<+\infty\), then \(M_T^Z\)
converges and so $M_t^Z/A_T\stackrel{a.s.}\to 0$. Otherwise, the Martingale Strong Law of Large Numbers 
\cite[Theorem~2.3.4, p.~46]{lapeyre2024montecarlo} gives
\(
{M_T^Z}/{(1+\langle M^Z\rangle_T)}\stackrel{a.s.}\longrightarrow 0.
\)
By the preceding estimate, on
\(\{A_T\to+\infty\}\) we have
\(
(1+\langle M^Z\rangle_T)/A_T=O(1).
\)
Hence, if \(\langle M^Z\rangle_\infty=+\infty\),
\[
\frac{M_T^Z}{A_T}
=
\frac{M_T^Z}{1+\langle M^Z\rangle_T}
\frac{1+\langle M^Z\rangle_T}{A_T}
\stackrel{a.s.}\longrightarrow 0\,.
\]
Thus \(M_T^Z=o(A_T)\) in either case.
The decomposition of \(Z_T\) gives
\(Z_T/A_T\stackrel{a.s.}\to 1\), proving
\ref{lem:compensated-increments:ii}.
For \ref{lem:compensated-increments:iii}, work on the event  \(\{Z_T\to+\infty\}\cap 
\{\E(\Delta Z_{t+1}\mid\F_t)/Z_t\to 0\}\). Part \ref{lem:compensated-increments:i} implies that
\(A_T\stackrel{a.s.}\to+\infty\). Set \(B_t:=1+A_t\). By
\ref{lem:compensated-increments:ii}, \(Z_t\stackrel{a.s.}\sim B_t\), and hence we have  
\(\E(\Delta Z_{t+1}\mid\F_t)/B_t\stackrel{a.s.}\to 0\). 
Consequently, 
since \(\E(\Delta Z_{t+1}\mid\F_t)\geq 0\) eventually with probability one, 
we have
\(
0\leq\E(\Delta Z_{t+1}\mid\F_t)\leq B_t
\)
for all sufficiently large \(t\), almost surely.
In particular, \(B_t>0\) and
\(
B_t\leq B_{t+1}
=B_t+\E(\Delta Z_{t+1}\mid\F_t)
\leq 2B_t
\)
eventually.
Then
\[
\frac{\E(\Delta Z_{t+1}\mid\F_t)}{B_t^2}
\leq
2\frac{\E(\Delta Z_{t+1}\mid\F_t)}{B_t(B_t+\E(\Delta Z_{t+1}\mid\F_t))}
=
2\left(\frac1{B_t}-\frac1{B_{t+1}}\right).
\]
It follows that \(\sum_t\E(\Delta Z_{t+1}\mid\F_t)/B_t^2<+\infty\) a.s.\,.
Because \(Z_t\stackrel{a.s.}\sim B_t\), this proves
\(\sum_t\E(\Delta Z_{t+1}\mid\F_t)/Z_t^2<+\infty\) a.s.\,. The conditional-variance
bound~\eqref{eq:lemma-cond-var-bound} gives the convergence of the second series  in \ref{lem:compensated-increments:iii}.
\end{proof}


\providecommand{\bysame}{\leavevmode\hbox to3em{\hrulefill}\thinspace}
\providecommand{\MR}{\relax\ifhmode\unskip\space\fi MR }
\providecommand{\MRhref}[2]{%
  \href{http://www.ams.org/mathscinet-getitem?mr=#1}{#2}
}
\providecommand{\href}[2]{#2}

\end{document}